\documentclass[10pt,twocolumn]{IEEEtran}

\IEEEoverridecommandlockouts                              

\usepackage{amsmath,graphicx,amsfonts,amssymb,amsthm,epsfig,mathrsfs,balance}
\usepackage{subcaption,tikz,color,multirow,algorithm,algpseudocode,algorithmicx}
\usepackage{cite,url,framed,bm,dsfont}

\usepackage{enumitem}

\newcommand{\spt}{{\mathrm{spt}}}

\newcommand{\differential}{{\rm{d}}}
\newcommand{\lr}{\ell_{\rm{rear}}}
\newcommand{\lf}{\ell_{\rm{front}}}

\let\emptyset\varnothing

\title{\LARGE \bf
Density-based Stochastic Reachability Computation for\\
Occupancy Prediction in Automated Driving
}

\author{Shadi Haddad, Abhishek Halder, and Baljeet Singh
\thanks{Shadi Haddad and Abhishek Halder are with the Department of Applied Mathematics, University of California, Santa Cruz, CA 95064, USA, {\texttt{\{shhaddad,ahalder\}@ucsc.edu}}}
\thanks{Baljeet Singh is with the Ford Greenfield Labs, Palo Alto, CA 94304, USA, {\tt\small{BSING124@ford.com}}}        
\thanks{This research was partially supported by 2019 Ford University Research Project, and the Chancellor Fellowship from the University of California, Santa Cruz.}
}

\begin{document}

\maketitle
\thispagestyle{empty}
\pagestyle{empty}

\begin{abstract}
We propose a stochastic reachability computation framework for occupancy prediction in automated driving by directly solving the underlying transport partial differential equation governing the advection of the closed-loop joint density functions. The resulting nonparametric gridless computation is based on integration along the characteristic curves, and allows online computation of the time-varying collision probabilities. Numerical simulations highlight the scope of the proposed method.
\end{abstract}




\section{Introduction}\label{SecIntro}
Guaranteeing safety without sacrificing performance in automated driving, especially in a mixed traffic of autonomous, semi-autonomous and human-driven vehicles, requires real-time planning and decision making against time-varying uncertainties. From a decision making (ego) vehicle’s perspective, these uncertainties stem not only from the environment (states and intents of neighboring traffic, road geometry, weather and lighting conditions), but also from the intrinsic sources (sensing and localization errors, communication packet drop, computational latency, model mismatch between the planning and control software). Thus, timely and accurate forecasting of the states of the ego vehicle as well as its environment is critical for provably safe decision making. At the same time, given the complexity of the dynamics and environment, it is extremely challenging to perform such forecasting computation at a time scale much smaller than the physical dynamics time scale. Therefore, a critical gap exists in real-time forecasting to enable provably safe automation for achieving objectives such as collision avoidance, safe lane change, and separation management.  This paper presents a framework for automated navigation by designing novel computational tools for stochastic uncertainty propagation enabling fast occupancy prediction in multi-lane driving scenarios. 

The importance of forecasts in automated driving while accounting the dynamic uncertainties has been well-advocated \cite{carvalho2015automated,hubmann2018automated} in the literature. The problem of stochastic occupancy prediction has been addressed before using finite state Markov chain abstractions \cite{althoff2009model,althoff2011comparison} on a discretization of the joint state and action space, and also by Monte Carlo methods \cite{broadhurst2005monte,eidehall2008statistical,funfgeld2017stochastic}. The purpose of the present work is to harness recent developments \cite{ehrendorfer1994liouville,halder2010beyond,halder2011dispersion,halder2015jgcd} in the theory and algorithms for solving the transport partial differential equation (PDE) governing the evolution of the joint state probability density function (PDF) subject to a known trajectory-level dynamics, to enable stochastic occupancy prediction in automated driving. This PDE-based computational framework has a number of attractive features pertinent to the occupancy prediction. For instance, the time-varying joint state PDFs can be computed \emph{exactly} along the trajectories of the closed-loop dynamics conditioned on the (joint) stochastic uncertainties in the initial conditions and parameters--a sharp juxtaposition with the discretization-based methods such as the Markov chain abstraction or the Monte Carlo, where the \emph{functional approximation} of the joint PDFs are computed--the quality and computational cost of approximation being dependent on the resolution of the discretization grid. In the proposed framework, the exact computation is possible since the underlying transport PDE is linear and its characteristic curves happen to coincide with the closed-loop trajectories, which can be algorithmically leveraged for scattered weighted point cloud-based computation, as we show in Section \ref{SecStocReachability}.

Another notable feature of the proposed framework is that the computation is nonparametric: we employ neither statistical (e.g., Gaussian, Gaussian mixture, exponential family) nor dynamical (e.g., polynomial vector field) approximation.  

The rest of this paper is organized as follows. 
In Section \ref{SecModels}, we describe two models for vehicle dynamics. Section \ref{SecStocReachability} explains the transport PDE-based stochastic reachability computation. Numerical simulations using the aforesaid models are detailed in Section \ref{SecNumericalSimulations}. Section \ref{secConclusions} concludes the paper.


\subsubsection*{Notations}
We use the bold-faced small letters for vectors, and bold-faced capital letters for matrices. The symbol ${\rm{diag}}$ denotes a digonal matrix, and ${\rm{blkdiag}}$ denotes a block-diagonal matrix. We use the notation $\bm{I}_{n}$ for the $n\times n$ identity matrix. For $\bm{x}\in\mathcal{X}\subseteq\mathbb{R}^{n}$, the notation $\bm{x} \sim \rho$ means that the random vector $\bm{x}$ has the joint PDF $\rho$. By definition,
\[\rho \geq 0\quad\text{for all}\quad\bm{x}\in\mathcal{X}, \quad\text{and}\quad\displaystyle\int_{\mathcal{X}}\rho\:\differential\bm{x} = 1.\]
The symbol $\mathcal{N}\left(\bm{\mu},\bm{\Sigma}\right)$ denotes a joint Gaussian PDF with the mean vector $\bm{\mu}$ and the covariance matrix $\bm{\Sigma}$. The support of a joint PDF $\rho$ is denoted as $\spt(\rho):=\{\bm{x}\mid\rho(\bm{x})\neq 0\}$. The notation $\mathds{1}_{\mathcal{X}}$ stands for the indicator function of a set $\mathcal{X}$. We use standard abbreviations such as ODE for ordinary differential equation, LTI for linear time invariant, MPC for model predictive control, QP for quadratic programming, and w.r.t. for ``with respect to".


\section{Models}\label{SecModels}
In this Section, we detail two vehicle dynamics models that will be used in this paper. The first is a kinematic bicycle model with four states and two controls, and the second is a dynamic bicycle model in road aligned coordinates with six states and three controls.

\begin{figure}[tpb]
        \centering
        \includegraphics[width=0.87\linewidth]{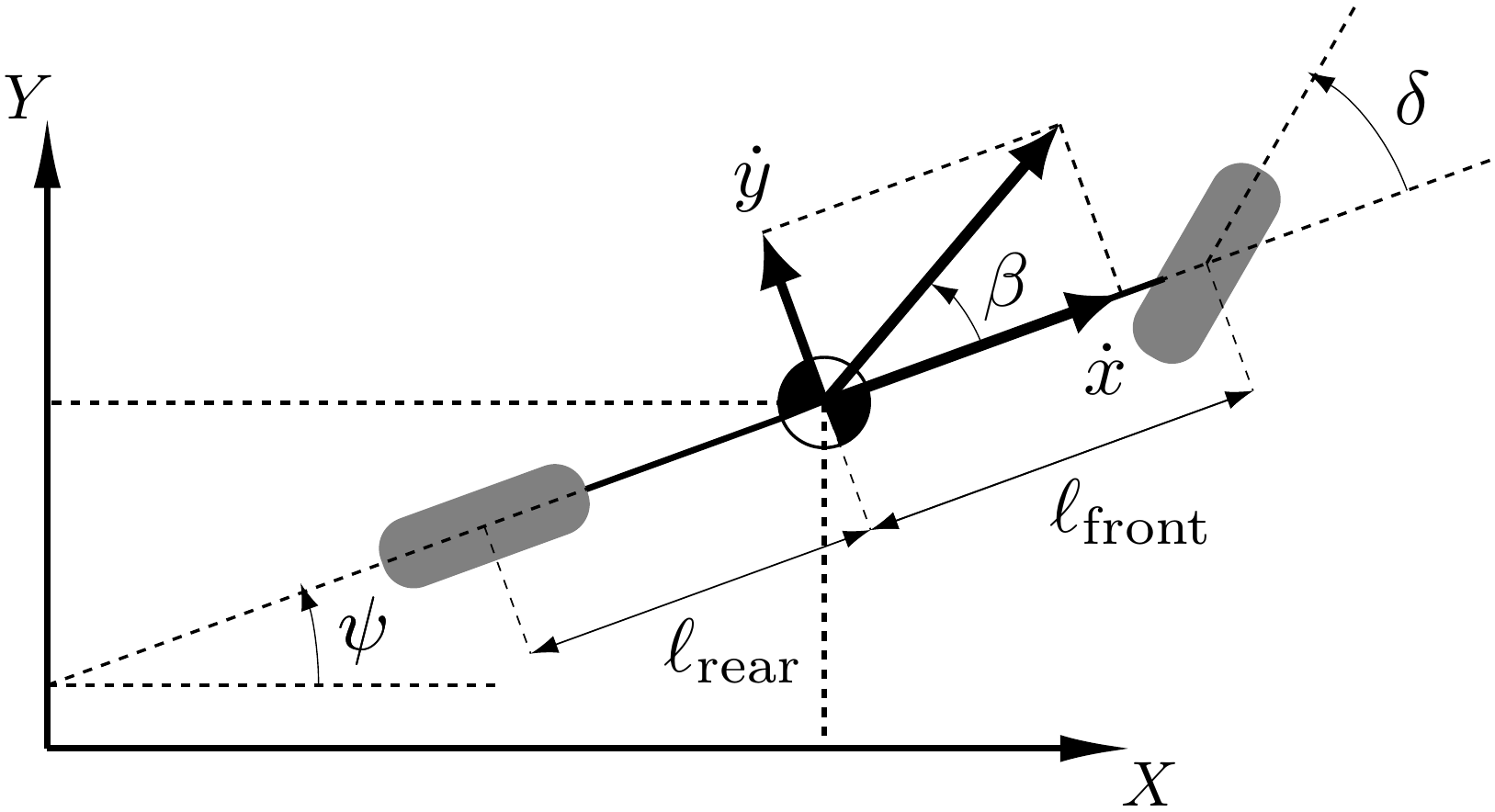}
        \caption{{\small{Schematic of notations used in the kinematic bicycle model in Section \ref{SubsecKinBicycleModel}. Here, $(X,Y)$ denotes an inertial coordinate system.}}}
\vspace*{-0.06in}
\label{Kinematicbicycle}
\end{figure}

\subsection{Kinematic Bicycle Model}\label{SubsecKinBicycleModel}
We consider the kinematic bicycle model \cite{rajamani2011vehicle} with state vector $\bm{x}:=\left(x,y,v,\psi\right)^{\!\top}$ and control vector $\bm{u}:=\left(a_{c},\delta\right)^{\!\top}$, given by (see Fig. \ref{Kinematicbicycle})
\begin{subequations}
\begin{align}
\dot{x} &= v\cos\left(\psi + \beta\right), \label{xdot}\\
\dot{y} &= v\sin\left(\psi + \beta\right), \label{ydot}\\
\dot{v} &= a_{c}, \label{vdot}\\
\dot{\psi} &= \frac{v}{\lr}\sin\beta, \label{psidot}
\end{align}
\label{KinematicBicyleModel}
\end{subequations}
wherein the sideslip angle
\begin{align}
\beta = \arctan\left(\frac{\lr}{\lf + \lr}\tan\delta\right),
\label{defbeta}
\end{align}
the parameters $\lf,\lr$ are the distances of the vehicle's center-of-mass to the front and rear axles, respectively. The state vector $\bm{x}$ comprises of the inertial position $(x,y)$ for the vehicle's center-of-mass, its speed $v$, and the vehicle's inertial heading angle $\psi$. The control vector $\bm{u}$ comprises of the acceleration $a_{c}$, and the front steering wheel angle $\delta$.

\begin{figure}[t]
        \centering
        \includegraphics[width=0.98\linewidth]{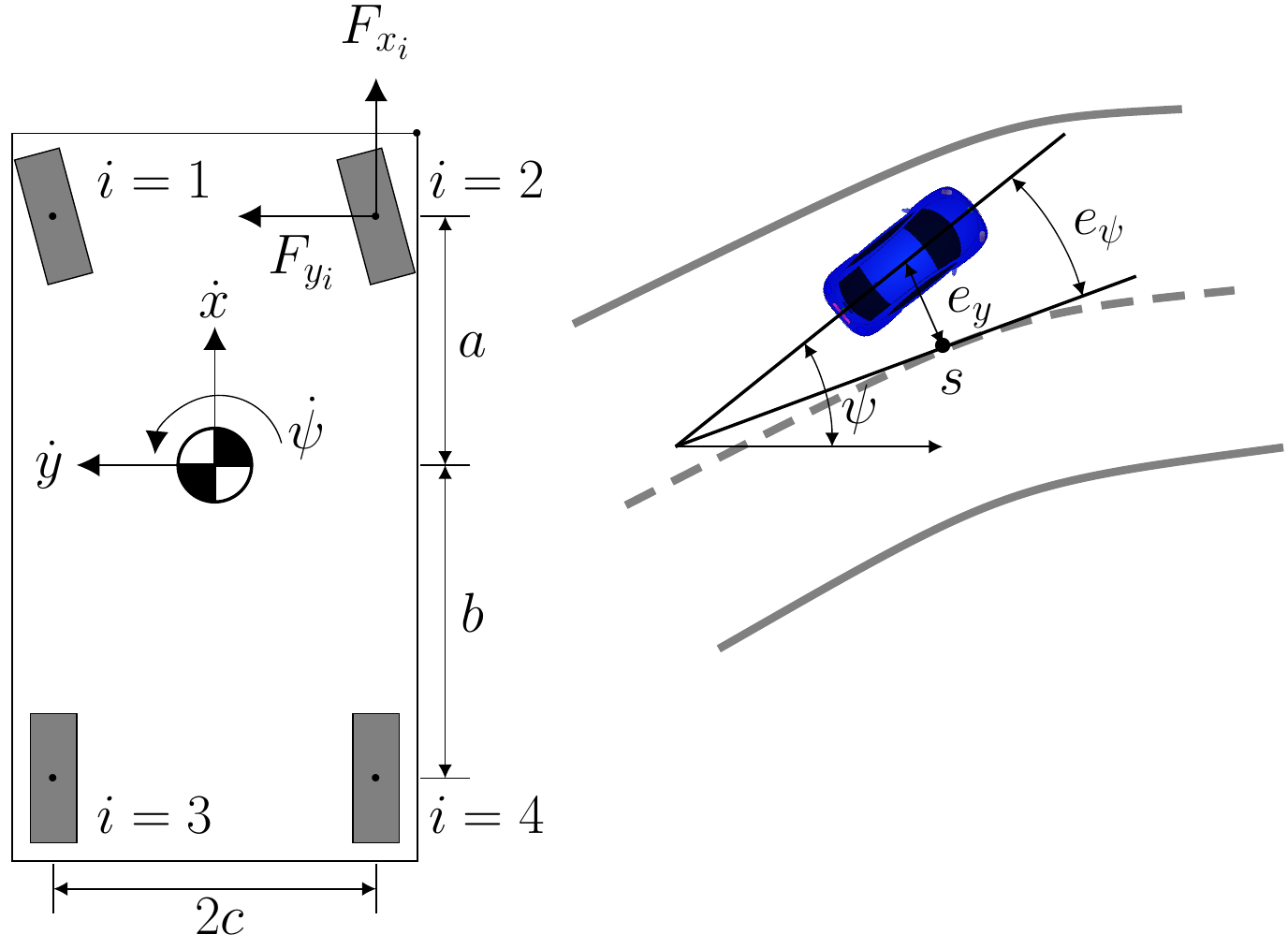}
        \caption{{\small{For the dynamic bicycle model in Section \ref{SubsecDynBicycleModel}, shown above are \emph{(left:)} the forces in the body frame, and \emph{(right:)} the vehicle states depicted in the road coordinates.}}}
\vspace*{-0.1in}
\label{States}
\end{figure}

\subsection{Dynamic Bicycle Model in Road Aligned Coordinates}\label{SubsecDynBicycleModel}
\begin{center}
\begin{table*}[h]
\centering
\begin{tabular}{| c | c | l |}
\hline
\:Symbol\: & \:Description\: & \:Numerical value\:\\
  \hline\hline
  & & \\			
  \:$a$\: & \: distance of the vehicle's center-of-gravity from its front axle \: & \:$1.432$ m\: \\
  & & \\
  \hline
  & & \\
  \:$b$\: & \: distance of the vehicle's center-of-gravity from its rear axle \: & \:$1.472$ m\: \\
  & & \\
  \hline
  & & \\
  \:$c$\: & \: the vehicle's half-track width \: & \:$0.8125$ m\: \\
  & & \\
  \hline 
   & & \\
  \:$m$\: & \: mass of the vehicle \: & \:$2050$ kg\: \\
  & & \\
  \hline
  & & \\
  \:$I_{z}$\: & \: yawing moment of inertia of the vehicle \: & \:$3344$ kg m$^2$\: \\
  & & \\
  \hline
\end{tabular}
\caption{The parameters used in the model (\ref{DynamicBicyleModelRoadAligned}).}
\label{ParamTable}
\end{table*}
\end{center}
\subsubsection{Open-loop dynamics}
We now consider the open-loop dynamic bicycle model \cite{Carvalho2013vehicle,kong2015kinematic}  in road aligned coordinates. The state vector is
\[\bm{x} := \left(v_{x},v_{y},v_{\psi},e_{\psi},e_{y},s\right)^{\!\top},\] 
where $v_{x}, v_{y}$ denote the longitudinal and lateral speeds, respectively, in m/s. Furthermore, $v_{\psi}$ denotes the yaw rate in rad/s. The states $\left(e_{\psi},e_{y},s\right)$ respectively denote the heading angle error in rad, the lateral position error measured w.r.t. the center of the road in m, and the longitudinal position of the vehicle along the road in m, as illustrated in Fig. \ref{States}. The control vector 
\[\bm{u}:=\left(\delta_{{\rm{front}}},\beta_{{\rm{left}}},\beta_{{\rm{right}}}\right)^{\!\top}\]
comprises of the front wheel steering angle $\delta_{{\rm{front}}}$, and the braking ratios for the left and right wheels: $\beta_{{\rm{left}}}$ and $\beta_{{\rm{right}}}$. The dynamics is given by
\begin{subequations}
\begin{align}
\dot{v}_{x} &= v_{y} v_{\psi} + \frac{1}{m}\displaystyle\sum_{i=1}^{4} F_{x_{i}}, \label{Dyn_vxdot}\\
\dot{v}_{y} &= -v_{x} v_{\psi} + \frac{1}{m}\displaystyle\sum_{i=1}^{4} F_{y_{i}},  \label{Dyn_vydot}\\
\dot{v}_{\psi} &= \frac{1}{I_{z}}\left(\!a\left(F_{y_{1}} + F_{y_{2}}\right) - b\left(F_{y_{3}} + F_{y_{4}}\right) + c\displaystyle\sum_{i=1}^{4} (-1)^{i}F_{x_{i}}\!\right), \label{Dyn_vpsidot}\\
\dot{e}_{\psi} &= v_{\psi} - \frac{\kappa}{1-\kappa e_{y}}\left(v_{x}\cos(e_{\psi}) - v_{y}\sin(e_{\psi})\right), \label{Dyn_epsidot}\\
\dot{e}_{y} &= v_{x}\sin(e_{\psi}) + v_{y}\cos(e_{\psi}),  \label{Dyn_eydot}\\
\dot{s} &= \frac{1}{1-\kappa e_{y}}\left(v_{x}\cos(e_{\psi}) - v_{y}\sin(e_{\psi})\right), \label{Dyn_sdot} 
\end{align}
\label{DynamicBicyleModelRoadAligned}
\end{subequations}


wherein the parameters pertaining to the vehicle geometry are listed in Table \ref{ParamTable}. Here, $\kappa$ denotes the curvature of the road.

\begin{figure}[t]
        \centering
        \includegraphics[width=0.98\linewidth]{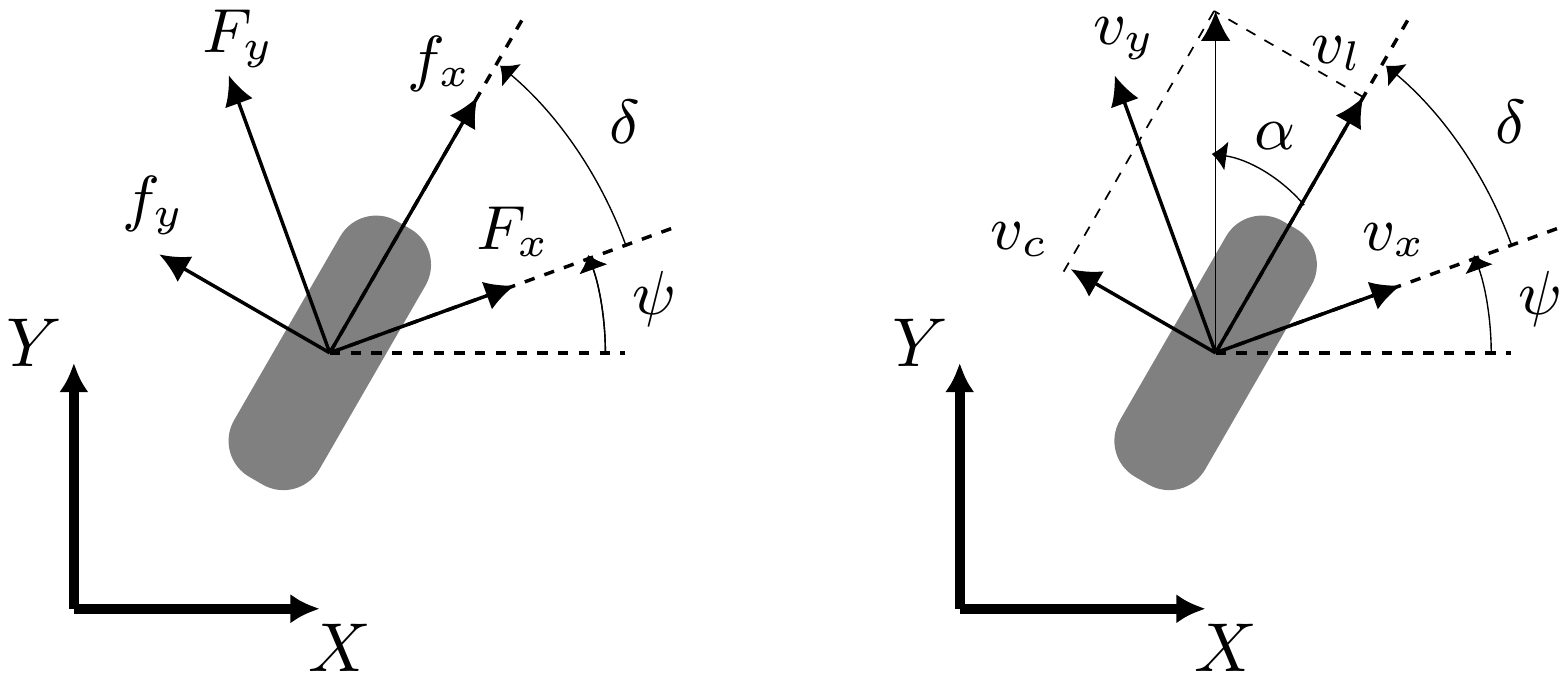}
        \caption{{\small{Schematic of the notations used in the dynamic bicycle model in Section \ref{SubsecDynBicycleModel}.}}}
\vspace*{-0.1in}
\label{Dynamicbicycle}
\end{figure}

To describe the external forces in (\ref{DynamicBicyleModelRoadAligned}), define the longitudinal and lateral forces (see Fig. \ref{Dynamicbicycle}), denoted by $f_{x_{i}}$ and $f_{y_{i}}$ respectively, as

\begin{subequations}
\begin{align}
f_{x_{i}} &:= \zeta \beta_{i} F_{z_{i}}, \label{fxi}\\
f_{y_{i}} &:= -C_{\alpha} \tan\alpha_{i}, \label{fyi}
\end{align}
\label{fxifyi}
\end{subequations}
where the braking ratio $\beta_{i} := \beta_{{\rm{left}}}$ for $i\in\{1,3\}$, and $\beta_{i} := \beta_{{\rm{right}}}$ for $i\in\{2,4\}$. In particular, $\beta_{i}\in[-1, 1]$ where $\beta_{i}=-1$ corresponds to maximum available braking, and $\beta_{i}=+1$ corresponds to maximum available throttle. We assume that the friction coefficient $\zeta$ between the tires and the road surface is the same for all the four tires, and remains constant over the time horizon of interest. In (\ref{fxi}), the normal forces $F_{z_{i}}$  are determined from the force and moment balance as
\begin{align}
F_{z_{1}} = F_{z_{2}} = \frac{mg}{2}\frac{b}{a+b}, \quad F_{z_{3}} = F_{z_{4}} = \frac{mg}{2}\frac{a}{a+b},
\label{Fzi}
\end{align}
where $g=9.81$ m/s$^{2}$ is the acceleration due to gravity. In (\ref{fyi}), the tire cornering stiffness parameter $C_{\alpha} = 250$ kN/rad, and the slip angles $\alpha_{i}$ satisfy
\begin{align} 
\tan\alpha_{i} = \frac{v_{c_{i}}}{v_{\ell_{i}}}, \quad i\in\{1,2,3,4\},
\label{tanalphai}
\end{align}

wherein $v_{c_{i}}$ and $v_{\ell_{i}}$, for $i\in\{1,2,3,4\}$, respectively denote the components of the $4\times 1$ vectors $\bm{v}_{c}$ and $\bm{v}_{\ell}$, given via the augmented equation
\begin{align}
\begin{pmatrix}
\bm{v}_{\ell}\\
\bm{v}_{c}
\end{pmatrix} = {\rm{blkdiag}}\left(\bm{U}(\delta_{{\rm{front}}}), \bm{U}(\delta_{{\rm{front}}}), \bm{I}_{2}, \bm{I}_{2}\right) \nonumber\\
\begin{pmatrix}
1 & 0 & -c\\
0 & 1 & a\\
1 & 0 & c\\
0 & 1 & a\\
1 & 0 & -c\\
0 & 1 &-b\\
1 & 0 & c\\
0 & 1 &-b
\end{pmatrix}
\begin{pmatrix}
v_{x}\\
v_{y}\\
v_{\psi}
\end{pmatrix}.
\label{vlvcvec}
\end{align} 
In (\ref{vlvcvec}), the orthogonal matrix
\begin{align}
\bm{U}(\delta_{{\rm{front}}}) := \begin{pmatrix}
\cos\left(\delta_{{\rm{front}}}\right) & \sin\left(\delta_{{\rm{front}}}\right)\\
-\sin\left(\delta_{{\rm{front}}}\right) & \cos\left(\delta_{{\rm{front}}}\right)
\end{pmatrix},
\label{DefMatUdeltafront}
\end{align}
where $\delta_{{\rm{front}}}$ is the first component of the control $\bm{u}$, and satisfies $-10\:\text{degrees} \leq \delta_{{\rm{front}}} \leq 10\:\text{degrees}$. 
The longitudinal and lateral forces, $F_{x_i}$ and $F_{y_i}$, respectively are defined in the body frame by
\begin{align}
\begin{pmatrix}
F_{x_{i}}\\
F_{y_{i}}
\end{pmatrix} = \bm{M}_{i} \begin{pmatrix}
f_{x_{i}}\\
f_{y_{i}}
\end{pmatrix},
\label{LongLateralForces}
\end{align}
where $\bm{M}_{i} = \bm{U}^{\top}\left(\delta_{{\rm{front}}}\right)$ for $i=1,2$, and $\bm{M}_{i} = \bm{I}_{2}$ for $i=3,4$. The forces $f_{x_{i}}$ and $f_{y_{i}}$ in (\ref{LongLateralForces}) are given by (\ref{fxifyi}).

\subsubsection{Feedback synthesis} \label{FeedbackSynthesis}
For feedback design, we first linearize the dynamics (\ref{DynamicBicyleModelRoadAligned}) around the ``velocity-hold" trim condition $v_{x} = 20$ m/s. We use \texttt{findop} in SIMULINK to compute the trim condition $(\bm{x}_{\text{trim}},\bm{u}_{\text{trim}})$ via sequential QP, where in addition to the equality constraint $v_{x} = 20$ m/s, we impose the inequality constraints:
\begin{subequations}
\begin{align}
\begin{pmatrix}
-10\:\text{degrees}\\
-1\\
-1
\end{pmatrix} \leq &\bm{u} \leq \begin{pmatrix}
10\:\text{degrees}\\
1\\
1
\end{pmatrix}, \label{ControlBounds}\\
\left(\mu_{e_{y}}\right)_{0} - 1.5 \leq & e_{y} \leq \left(\mu_{e_{y}}\right)_{0} + 1.5. \label{eyBounds}
\end{align}
\label{TrimConstraints}
\end{subequations}
The constraints (\ref{ControlBounds}) enforce bounded controls. The path constraint (\ref{eyBounds}) enforces the lateral position error $e_{y}$ to be within $\pm 1.5$ m of the initial mean lateral position error $\left(\mu_{e_{y}}\right)_{0}$. Since we allow stochastic uncertainties in initial condition, the mean, and hence the computation of the trim will depend on the joint state PDF at time $t=0$.

Linearizing (\ref{DynamicBicyleModelRoadAligned}) around the resulting trim point $(\bm{x}_{\text{trim}},\bm{u}_{\text{trim}})$, we obtain the LTI matrix pair $\left(\bm{A}_{\text{trim}},\bm{B}_{\text{trim}}\right)$, which is then used to obtain the explicit MPC feedback $\bm{u} = \bm{\pi}(\bm{x},t)$ for this linearized system to minimize the deviation from $(\bm{x}_{\text{trim}},\bm{u}_{\text{trim}})$ over a prediction horizon subject to the constraints (\ref{TrimConstraints}) via the MPC toolbox \cite{herceg2013multi}. Specifically, let $\Delta\bm{x} := \bm{x}- \bm{x}_{\text{trim}}$, $\Delta\bm{u} := \bm{u} - \bm{u}_{\text{trim}}$. The MPC objective penalizes the cumulative cost

\begin{equation}
\begin{split}
&\!\!\!\int_{t}^{t+t_{p}}\!\!\!\left(\Delta\bm{x}(\tau)\right)^{\top}\bm{Q}\Delta\bm{x}(\tau) \:\differential\tau\\
 &\!\!\!\!\!+\! \int_{t}^{t+t_{p}}\!\!\!\left(\left(\Delta\bm{u}(\tau)\right)^{\top}\bm{R}\Delta\bm{u}(\tau)+ \left(\Delta\dot{\bm{u}}(\tau)\right)^{\top}\bm{S}\Delta\dot{\bm{u}}(\tau)\right)\:\differential\tau,
\end{split}
\label{MPCobjective}
\end{equation}
over a prediction horizon $t_{p}$ subject to the trimmed LTI dynamics and the aforesaid constraints; the weight matrices $\bm{Q},\bm{R},\bm{S}$ correspond to state penalty, control penalty, and slew rate penalty, respectively.

As is well-known \cite{bemporad2002explicit}, computing the explicit MPC feedback for an LTI system amounts to solving parametric QP, whose solution is a continuous piecewise affine policy
\begin{align}
 \bm{\pi}(\bm{x},t) = \bm{\Gamma}_{j}\bm{x} + \bm{\gamma}_{j}, \quad\text{for}\quad \bm{x}(t)\in\mathcal{P}_{j},
\label{ExplicitMPCfeedback}
\end{align}
where $j=1,2,\hdots,\nu$, and the disjoint polytopic partition $\displaystyle\cup_{j=1}^{\nu}\mathcal{P}_{j}$ of the reach set of the closed-loop constrained LTI system depends on the time horizon in MPC synthesis, as well as on the constraints (\ref{TrimConstraints}) and the pair $\left(\bm{A}_{\text{trim}},\bm{B}_{\text{trim}}\right)$. We perform the explicit MPC synthesis offline, and store the feedback gain matrix-vector pairs $\{\bm{\Gamma}_{j},\bm{\gamma}_{j}\}_{j=1}^{\nu}$ as well as the polytopic partition $\displaystyle\cup_{j=1}^{\nu}\mathcal{P}_{j}$.

In Section \ref{SimulationDynamic}, we will discuss a multi-lane driving scenario where different vehicles will have different initial joint state PDFs, thus different trim points, and consequently different MPC feedback policies.

\begin{figure}[t]
        \centering
        \includegraphics[width=0.98\linewidth]{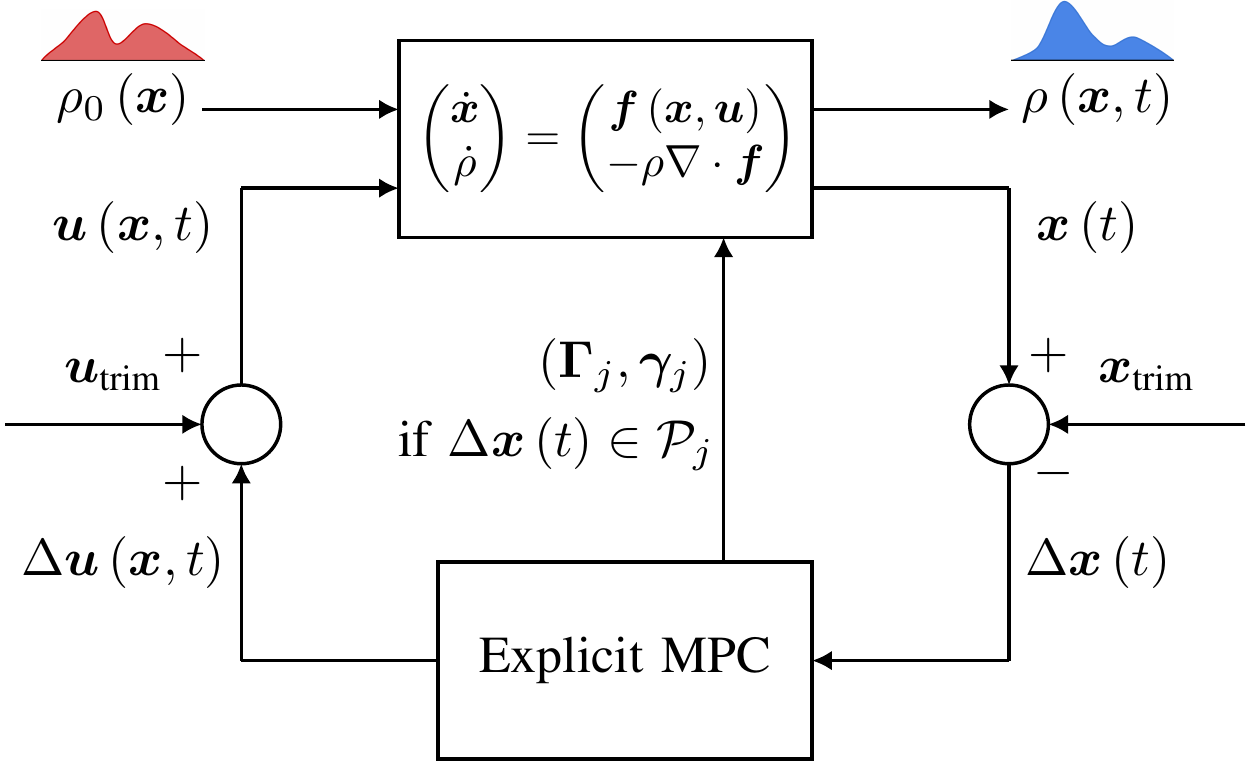}
        \caption{{\small{Block diagram for the computation of the time-varying joint state PDF $\rho(\bm{x},t)$ with open-loop dynamics $\dot{\bm{x}}=\bm{f}(\bm{x},\bm{u})$ and explicit MPC feedback, subject to the initial joint state PDF $\rho_{0}(\bm{x})$. The linearized explicit MPC is synthesized about the trim state-control pair $(\bm{x}_{\text{trim}},\bm{u}_{\text{trim}})$ (see Section \ref{FeedbackSynthesis}).}}}
\vspace*{-0.15in}
\label{BlockDiagram}
\end{figure}

\section{Stochastic Reachability via Liouville PDE}\label{SecStocReachability}
The control system (\ref{KinematicBicyleModel})-(\ref{defbeta}) is of the standard form $\dot{\bm{x}} = \bm{f}\left(\bm{x},\bm{u},t\right)$ with initial condition $\bm{x}_{0}:=\bm{x}(t=0)$. We would like to compute the evolution of the transient joint state PDF $\rho(\bm{x},t)$ subject to this controlled dynamics for a given non-randomized control policy $\bm{u} = \bm{\pi}\left(\bm{x},t\right)$, and uncertainties in the initial conditions $\bm{x}_{0} \sim \rho_{0}$ (given). Denoting the closed-loop dynamics as

\begin{align}
\dot{\bm{x}} = \bm{f}\left(\bm{x},\bm{\pi}\left(\bm{x},t\right),t\right) =: \bm{g}\left(\bm{x},t\right), \quad \bm{x}(t=0) = \bm{x}_{0},
\label{ClosedLoopODE}
\end{align} 
the associated dynamics of the joint state PDF is governed by the Liouville PDE initial value problem (IVP)
\begin{align}
\dfrac{\partial\rho}{\partial t} = -\nabla\cdot\left(\rho\bm{g}\right), \quad \rho\left(\bm{x},t=0\right) = \rho_{0}\left(\bm{x}\right),
\label{LiouvillePDE}
\end{align}
where $\nabla\cdot$ denotes the divergence operator w.r.t. the state vector $\bm{x}$, and the function $\rho_{0}$ is the joint PDF describing the initial condition uncertainties. The solution of the IVP (\ref{LiouvillePDE}) is the transient joint state PDF $\rho(\bm{x},t)$ which characterizes the time varying state uncertainties subject to the closed-loop dynamics (\ref{ClosedLoopODE}). Our interest is to perform fast nonparametric computation of $\rho(\bm{x},t)$.

Since the IVP (\ref{LiouvillePDE}) involves a first order PDE,  it is possible to solve the same using the method-of-characteristics. Specifically, it can be shown \cite[Sec. II.A]{halder2011dispersion} that the characteristic curves associated with (\ref{LiouvillePDE}) are the flows generated by (\ref{ClosedLoopODE}), and hence must be non-intersecting due to standard existence-uniqueness of the flow for an ODE of the form (\ref{ClosedLoopODE}) with smooth vector field $\bm{g}$. This can be exploited to design a weighted scattered point cloud-based computation of the transient joint PDFs $\rho(\bm{x},t)$ by integrating the following system of characteristic ODEs (see \cite{halder2011dispersion}) w.r.t. time $t$:

\begin{align}
\begin{pmatrix}
\dot{\bm{x}}^{i}\\
\dot{\rho}^{i}
\end{pmatrix} \!= 
\begin{pmatrix}
\bm{g}\left(\bm{x}^{i},t\right)\\
-\rho^{i}\nabla\cdot\bm{g}\left(\bm{x}^{i},t\right)
\end{pmatrix}, \quad \begin{pmatrix}
\bm{x}^{i}(t=0)\\
\rho^{i}(t=0)
\end{pmatrix} \!= \begin{pmatrix}
\bm{x}_{0}^{i}\\
\rho_{0}^{i}
\end{pmatrix},
\label{CharODE}
\end{align}

where the superscript $i=1,\hdots,N$ denotes the sample index. In other words, one can first generate $N$ random samples from the known initial joint state PDF $\rho_{0}$, and then the ODE (\ref{CharODE}) can be integrated \emph{along the characteristics} for each of the $N$ samples, resulting in a weighted scattered point cloud representation $\{\bm{x}^{i}(t),\rho^{i}(t)\}_{i=1}^{N}$ of the transient joint state PDF $\rho(\bm{x},t)$. For $i=1,\hdots,N$, the vector $\bm{x}^{i}(t)$ tells the location of the $i$-th sample in the state space at time $t$, and the nonnegative weight $\rho^{i}(t)$ gives the value of the joint state PDF evaluated at that state space location at that time. Thus, at any given time $t>0$, the higher (resp. lower) value of $\rho^{i}$ quantifies the higher (resp. lower) likelihood of the state $\bm{x}^{i}$.

Several features of the scattered point cloud-based computation of the joint PDF $\rho(\bm{x},t)$ described above make it particularly attractive for real-time occupancy prediction in automated driving scenarios:

\begin{enumerate} [label=(\roman*)]
\item The computation along each characteristic curve is independent of the other, and hence the integration of (\ref{CharODE}) is massively parallelizable. 
\item The computation is nonparametric, and does not suffer from the moment closure problem. From a practical standpoint, one cannot guarantee a priori if the joint PDF evolution subject to (\ref{ClosedLoopODE}) may admit a finite-dimensional sufficient statistic, and hence a nonparametric computation is desired. From an information-theoretic viewpoint, nonparametric computation (as opposed to moment-based approximation) is preferred since the joint PDF subsumes all of its moment information. 
\item Unlike traditional Monte-Carlo methods, the proposed method explicitly computes $\rho^{i}(t)$ via the characteristic ODE in (\ref{CharODE}). This allows gridless computation. In contrast, traditional Monte-Carlo methods only propagate the states $\bm{x}^{i}(t)$, and then approximates the joint state PDF from the state samples alone. This function approximation requires discretizing the state space and suffers from the ``curse of dimensionality" \cite{bellmanBook}. The explicit time integration of the characteristic ODE circumvents this difficulty.
\end{enumerate}

We remark here that it is possible to use the characteristic ODE (\ref{CharODE}) to write the solution of the PDE IVP (\ref{LiouvillePDE}) in semi-analytical form:

\[\rho(\bm{x},t) = \rho_{0}(\bm{x}_{0}\left(\bm{x},t\right))\exp\left(-\int_{0}^{t}\left(\nabla\cdot\bm{g}\right)(\tau)\:\differential\tau\right),\]
where $\bm{x}_{0}\left(\bm{x},t\right)$ is the inverse flow map for the closed-loop dynamics (\ref{ClosedLoopODE}). In most vehicular models such as those described in Section \ref{SecModels}, one does not have analytical handle on this inverse flow map, and hence it is computationally beneficial to directly integrate the system of characteristic ODEs (\ref{CharODE}) along the sample trajectories. In the numerical simulations described next, we do so using the standard \texttt{ode45} integrator in MATLAB\textsuperscript{\textregistered}.


\section{Numerical Simulations}\label{SecNumericalSimulations}
In the following, we elucidate the overall computational framework by first giving some illustrative examples in Section \ref{subsecIllustrativeExamples}, and then consider a case study in Section \ref{subsecCaseStudy} to highlight the scope of the proposed method. Section \ref{subsecComparisonMC} provides a comparison of the numerical performance of the proposed method with the standard Monte Carlo. All simulations were performed in MATLAB\textsuperscript{\textregistered} R2019b running on iMac with 3.4 GHz Intel Core i5 processor and 8 GB memory. 

\subsection{Illustrative Examples}\label{subsecIllustrativeExamples}

We consider a 2-lane unidirectional highway driving scenario with 2 vehicles, one ego and another non-ego, driving in the adjacent lanes with stochastic uncertainties in their states. The ego vehicle has better estimation accuracy of its own states compared to the same for the state of the non-ego vehicle. Consequently, the uncertainties in the state of the non-ego vehicle at any time, as perceived by the ego vehicle, is expected to be more than the uncertainties in its own state at that time. We often refer to the ego vehicle as the ``vehicle 1", and the non-ego vehicle as the ``vehicle 2".

We suppose that the time horizon of interest is $[t_{0}, t_{f}]$. Supposing that at $t_{0}$, the ego vehicle has an estimate of its own as well as the non-ego vehicle's state in terms of the respective beliefs or joint state PDFs: $\rho_{0}^{\text{ego}}(\bm{x})$ and $\rho_{0}^{\text{non-ego}}(\bm{x})$, possibly computed via some sensor fusion algorithm using the particle filter or the extended Kalman filter, the ego vehicle would like to predict the transient beliefs $\rho^{\text{ego}}(\bm{x},t)$ and $\rho^{\text{non-ego}}(\bm{x},t)$ for $t\in[t_{0}, t_{f}]$. The horizon length $t_{f} - t_{0}$ is typically few seconds, and the ego vehicle may repeat this predictive computation once new estimate (i.e., belief) arrives at the end of this horizon. We only illustrate the computation for a single horizon, set $t_{0}\equiv 0$ without loss of generality, and refer to $\bm{x}(t_{0})$ as the ``initial" condition.

To compute $\rho^{\text{ego}}(\bm{x},t)$ and $\rho^{\text{non-ego}}(\bm{x},t)$ for $t\in[t_{0}, t_{f}]$, conditioned on the initial beliefs $\rho_{0}^{\text{ego}}(\bm{x})$ and $\rho_{0}^{\text{non-ego}}(\bm{x})$ at $t_{0}$, the ego vehicle utilizes its own as well as the non-ego vehicle's closed-loop dynamics. The latter may only be an approximation since the exact non-ego dynamics is not known to the ego vehicle, but the approximation may be refined and/or learnt from observations in real time. We illustrate the proposed framework using the models described in Section \ref{SecModels}.

\subsubsection{Stochastic reachability with the kinematic bicycle model}\label{SimulationKinematic}
Let us consider the case where both the ego and the non-ego vehicles have the dynamics as in Section \ref{SubsecKinBicycleModel} with $\ell_{\text{front}}=1$ m and $\ell_{\text{rear}}=1.5$ m. Suppose that at $t_{0}=0$, their initial joint state PDFs are
\begin{align}
\rho_{0}^{\text{ego}}(\bm{x}) = \mathcal{N}\left(\bm{\mu}_{1},\bm{\Sigma}_{1}\right), \quad \rho_{0}^{\text{non-ego}}(\bm{x}) = \mathcal{N}\left(\bm{\mu}_{2},\bm{\Sigma}_{2}\right),
\label{initialJointGaussian}	
\end{align}
with
\begin{align*}
\bm{\mu}_{1} &= \left(0, 0, 20, 0\right)^{\top}, \quad \bm{\Sigma}_{1} = {\rm{diag}}\left(10^{-2}, 10^{-2}, 10^{-1}, 10^{-3}\right),\\
\bm{\mu}_{2} &= \left(0, 5, 20, 0\right)^{\top},\quad \bm{\Sigma}_{2} = {\rm{diag}}\left(10^{-2}, 10^{-1}, 1, 10^{-1}\right).	
\end{align*}
For all $t\in [t_{0},t_{f}] \equiv [0,5]$, we set the acceleration input $a_{c}(t) \equiv \sin(t)$, and the front steering wheel angle input $\delta(t) \equiv 0$ for both the vehicles. We generate $N=1000$ random samples for each of the above two initial joint state PDFs, evaluate the resulting samples exactly at the respective known initial joint state PDFs, and then propagate the two weighted scattered point clouds $\{\bm{x}^{\text{ego},i}(t),\rho^{\text{ego},i}(t)\}_{i=1}^{N}$ and $\{\bm{x}^{\text{non-ego},i}(t),\rho^{\text{non-ego},i}(t)\}_{i=1}^{N}$ as described in Section \ref{SecStocReachability}. As explained before, this joint PDF propagation via the method-of-characteristics is a gridless computation.

\begin{figure*}[p]
        \centering
        \includegraphics[width=0.85\linewidth]{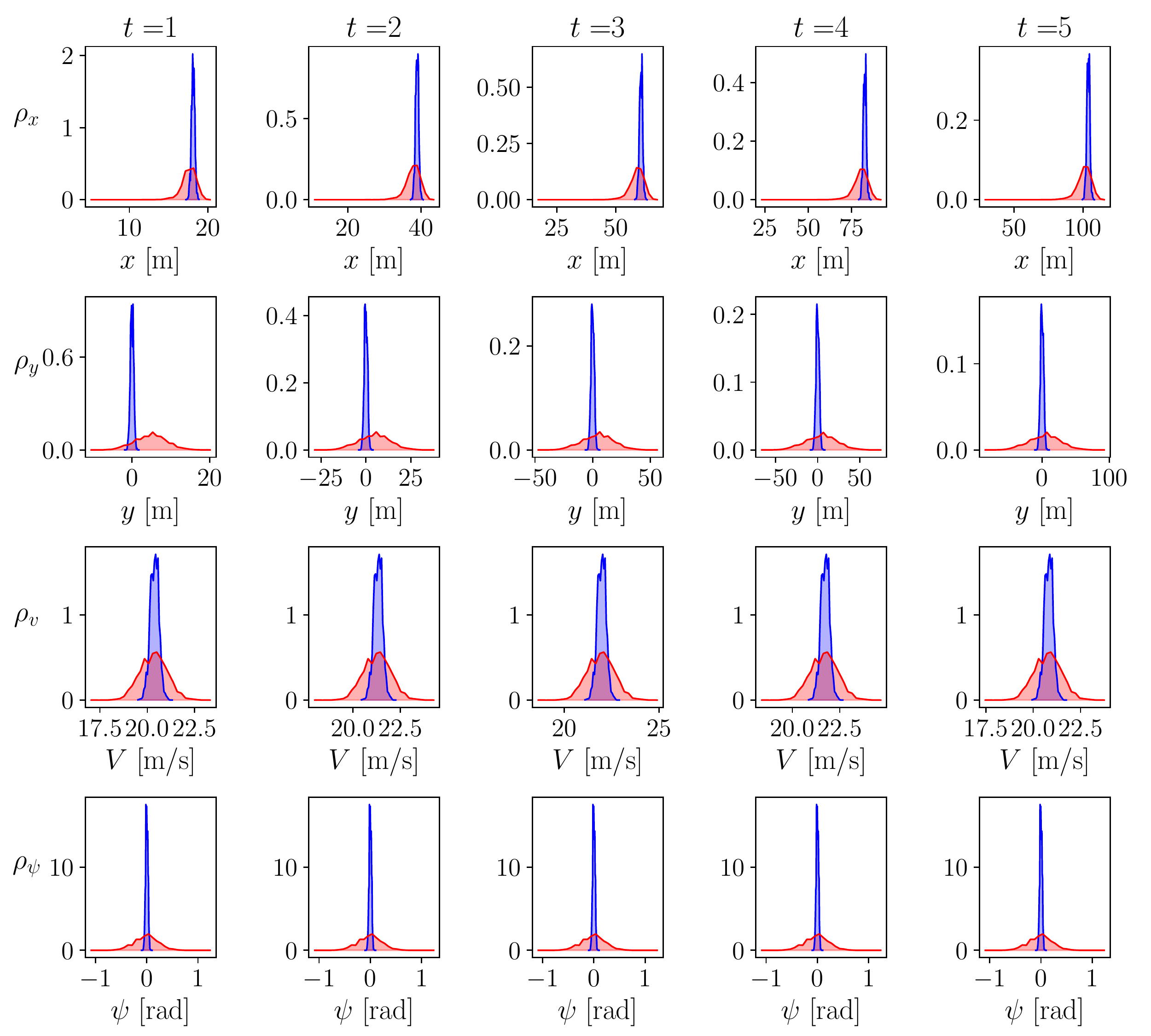}
               \caption{{\small{For the simulation set up detailed in Section \ref{SimulationKinematic}, the evolution of the transient univariate marginal state PDFs for the vehicle 1 (ego, \emph{in blue}) and vehicle 2 (non-ego, \emph{in red}).}}}
\vspace*{-0.05in}
\label{Kinematic1D}

\vspace*{.3 in}
        \includegraphics[width=0.85\linewidth]{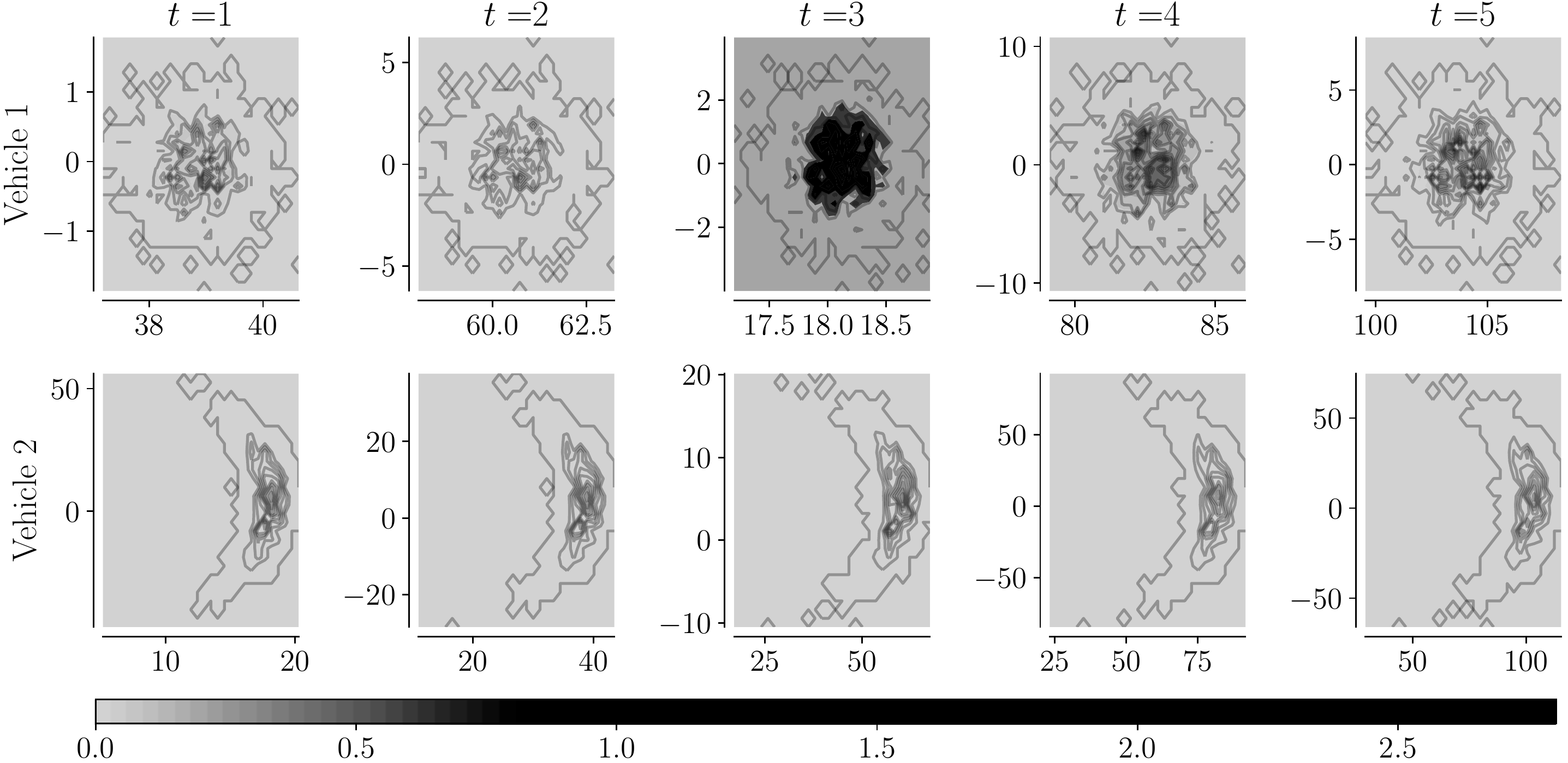}
         \caption{{\small{For the simulation set up detailed in Section \ref{SimulationKinematic}, the evolution of the transient bivariate marginal state PDF in $\left(x,y\right)$ variables for the vehicle 1 (ego) and vehicle 2 (non-ego). The colorbar (\emph{light hue = small}, \emph{dark hue = large}) shows the values of the bivariate marginals.}}}
\vspace*{-0.05in}
\label{KinematicJoint}

\end{figure*}

Fig. \ref{Kinematic1D} shows the evolution of the univariate marginal state PDFs for the two vehicles, obtained by numerically integrating (see e.g., \cite[Sec. V.B.3]{halder2011dispersion}) the weighted scattered point cloud data described above. The evolution of the bivariate marginal PDFs in $\left(x,y\right)$ variables is depicted in Fig. \ref{KinematicJoint}. These plots indeed show that the uncertainties in the (non-ego) vehicle 2's state  at any given time, has more dispersion than the (ego) vehicle 1's state.


\begin{figure*}[p]
        \centering
        \includegraphics[width=0.95\linewidth]{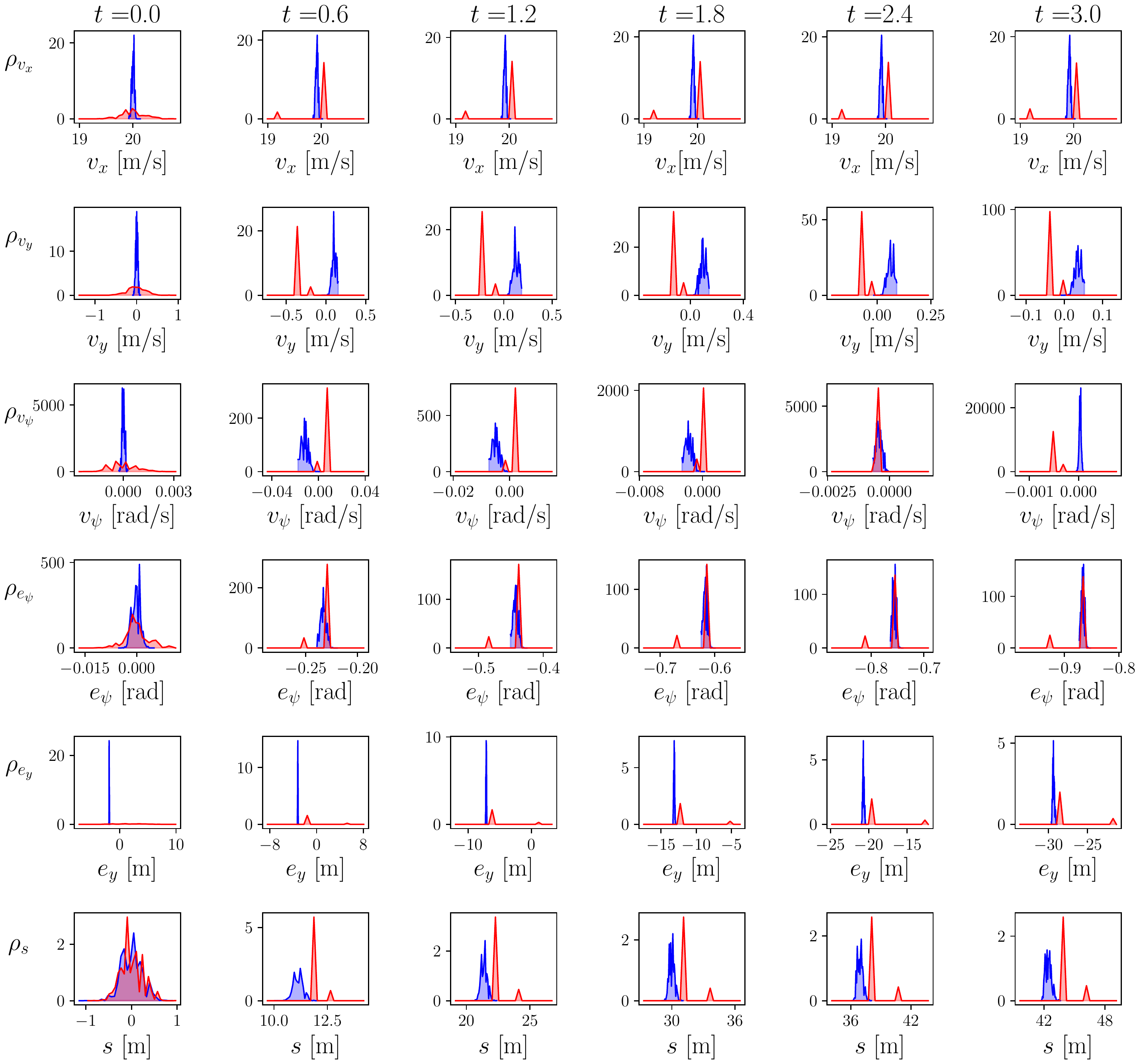}
         \caption{{\small{For the simulation set up detailed in Section \ref{SimulationDynamic}, the evolution of the transient univariate marginal state PDFs for the vehicle 1 (ego, \emph{in blue}) and vehicle 2 (non-ego, \emph{in red}).}}}
\vspace*{-0.05in}
\label{Dynamic1D}

\vspace*{0.1 in}

        \includegraphics[width=0.95\linewidth]{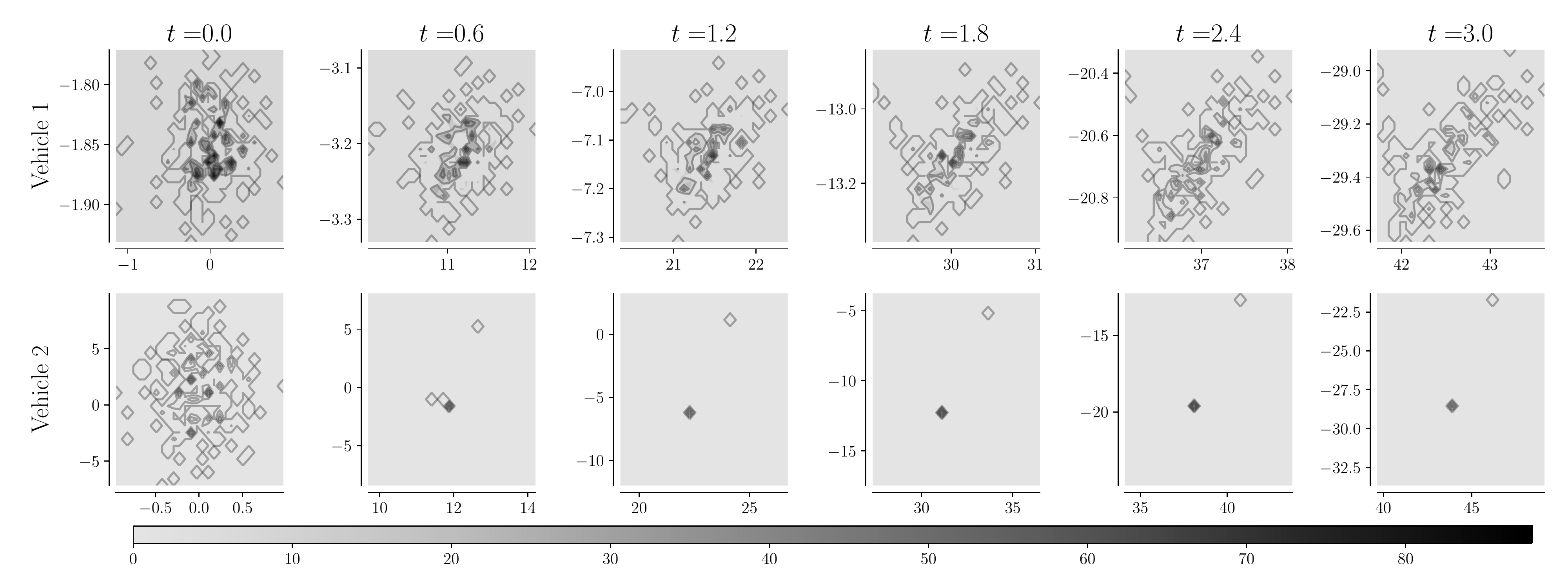}
         \caption{{\small{For the simulation set up detailed in Section \ref{SimulationDynamic}, the evolution of the transient bivariate marginal state PDF in $\left(s,e_y\right)$ variables for the vehicle 1 (ego) and vehicle 2 (non-ego). The colorbar (\emph{light hue = small}, \emph{dark hue = large}) shows the values of the bivariate marginals.}}}
\vspace*{-0.05in}
\label{DynamicJoint}
\end{figure*}


\subsubsection{Stochastic reachability with the dynamic bicycle model}\label{SimulationDynamic}
We now consider the case where both the ego and the non-ego vehicles have the dynamics as in Section \ref{SubsecDynBicycleModel} with trim-linearized explicit MPC feedback satisfying the constraints (\ref{TrimConstraints}). As we mentioned in Section \ref{FeedbackSynthesis}, the computation of the trim depends on the initial joint state PDF due to the constraint (\ref{eyBounds}). We suppose that the initial joint state PDFs are of the form (\ref{initialJointGaussian}) with 
\begin{align*}
\bm{\mu}_{1} &= \left(20, 0, 0, 0, -1.85, 0\right)^{\top}, \\
\bm{\Sigma}_{1} &= {\rm{diag}}\left(0.11, 0.11, 1.24\times10^{-8}, 2.78\times10^{-6}, 10^{-2}, 0.11\right),\\
\bm{\mu}_{2} &= \left(20, 0, 0, 0, 1.85, 0\right)^{\top},\\
\bm{\Sigma}_{2} &= {\rm{diag}}\left(1.11, 1.11, 1.24\times10^{-7}, 2.5\times10^{-5}, 1.11, 0.11\right).	
\end{align*}

We perform the offline linearized explicit MPC synthesis using the MPC toolbox \cite{herceg2013multi} with sampling time 0.1 s, prediction horizon $t_{p} = 3$ s, and control horizon of length 2 s. The weight matrices in (\ref{MPCobjective}) are chosen as below so that the vehicle closed loop response for speed control and lane centering has the desired settling time of $O(1)$ seconds:
\begin{align*}
\bm{Q} = 10\bm{I}_{6}, \quad \bm{R}= \bm{I}_{3}, \quad \bm{S} = 10^{-1}\bm{I}_{3}.	
\end{align*}
The resulting continuous piecewise affine feedback (\ref{ExplicitMPCfeedback}) admits a disjoint polytopic partition $\cup_{j=1}^{3610}\mathcal{P}_{j}$ in the six dimensional state space. In our simulation environment mentioned before, the offline MPC synthesis takes approximately 35 minutes.

We generate $N=200$ random samples for each of the two initial joint state PDFs, and as before, evaluate the resulting samples exactly at the respective known initial joint state PDFs. For $t\in [t_{0},t_{f}] \equiv [0,3]$, we propagate the two weighted scattered point clouds $\{\bm{x}^{\text{ego},i}(t),\rho^{\text{ego},i}(t)\}_{i=1}^{N}$ and $\{\bm{x}^{\text{non-ego},i}(t),\rho^{\text{non-ego},i}(t)\}_{i=1}^{N}$ for the closed-loop dynamics with respective explicit MPC feedback. Thus, the joint PDF computation depends on the feedback gain matrix-vector pairs $\{\bm{\Gamma}_{j},\bm{\gamma}_{j}\}_{j=1}^{\nu}$ and the polytopic partition $\cup_{j=1}^{\nu}\mathcal{P}_{j}$ (in our simulation, $\nu = 3610$).

Fig. \ref{Dynamic1D} shows the evolution of the univariate marginal state PDFs for the two vehicles, obtained by numerically integrating the weighted scattered point cloud data described above. The evolution of the bivariate marginal PDFs in $\left(s,e_y\right)$ variables is depicted in Fig. \ref{DynamicJoint}. From these plots, we observe that the uncertainties in the (non-ego) vehicle 2's state  at any given time, has more dispersion than the (ego) vehicle 1's state, as expected. We next use the bivariate $\left(s,e_y\right)$ marginal PDFs to compute the time-varying collision probability over the simulation interval $[t_{0},t_{f}] \equiv [0,3]$.


\begin{figure}[t]
        \centering
        \includegraphics[width=0.9\linewidth]{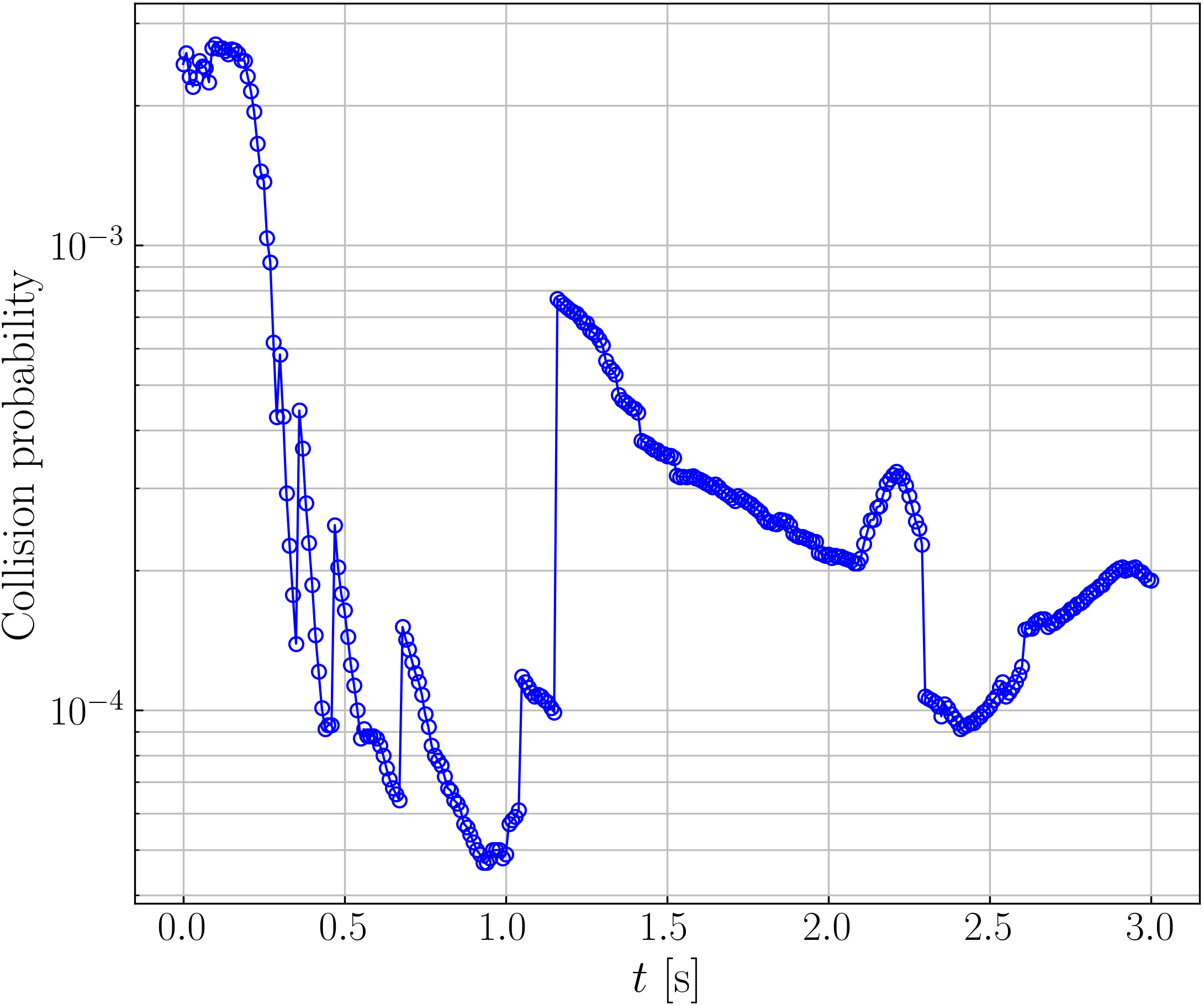}
         \caption{{\small{The collision probability $p_{\text{collision}}(t)$ between the ego and the non-ego vehicle for the simulation set up detailed in Section \ref{SimulationDynamic}.}}}
\vspace*{-0.1in}
\label{FigCollision}
\end{figure}


Notice that at any time $t\in[t_{0},t_{f}]$, nonzero collision probability results only if the supports of the $(s,e_y)$ bivariate marginal PDFs for the ego and the non-ego vehicles intersect at that time. Let us denote these bivariate marginal PDFs as $\xi^{\text{ego}}\left(s^{\text{ego}},e_{y}^{\text{ego}},t\right)$ and $\xi^{\text{non-ego}}\left(s^{\text{non-ego}},e_{y}^{\text{non-ego}},t\right)$, respectively. Making the standard assumption \cite[Sec. IV]{lambert2008collision} that the dynamics of the the ego and the non-ego vehicles are independent in $[t_{0},t_{f}]$, the collision probability $p_{\text{collision}}(t)$ is given by
\begin{align}
p_{\text{collision}}(t) = \!\displaystyle\int_{\mathcal{D}(t)}\!\xi^{\text{ego}}&\left(s^{\text{ego}},e_{y}^{\text{ego}},t\right)\xi^{\text{non-ego}}\left(s^{\text{non-ego}},e_{y}^{\text{non-ego}},t\right)\nonumber\\
&\differential s^{\text{ego}}\differential e_{y}^{\text{ego}}\differential s^{\text{non-ego}}\differential e_{y}^{\text{non-ego}}, 
\label{CollisionProb}	
\end{align}
where the time-varying domain
\begin{align*}
&\mathcal{D}(t):=\big\{\left(s^{\text{ego}},e_{y}^{\text{ego}},s^{\text{non-ego}},e_{y}^{\text{non-ego}}\right)\in\mathbb{R}^{4} \mid \\
&\spt\!\left(\xi^{\text{ego}}\!\left(s^{\text{ego}},e_{y}^{\text{ego}},t\right)\right)\!\bigcap\spt\!\left(\xi^{\text{non-ego}}\!\left(s^{\text{non-ego}},e_{y}^{\text{non-ego}},t\right)\right)\!\neq\emptyset \}.
\end{align*}
We can rewrite (\ref{CollisionProb}) as
\begin{align}
p_{\text{collision}}(t) = \!\displaystyle\int_{\mathbb{R}^{4}}\!\!\!\mathds{1}_{\mathcal{D}(t)}\xi^{\text{ego}}\!&\left(s^{\text{ego}},e_{y}^{\text{ego}},t\right)\!\xi^{\text{non-ego}}\!\left(s^{\text{non-ego}},e_{y}^{\text{non-ego}},t\right)\nonumber\\
&\differential s^{\text{ego}}\differential e_{y}^{\text{ego}}\differential s^{\text{non-ego}}\differential e_{y}^{\text{non-ego}}.
\label{CollisonProbIndicator}	
\end{align}
In \cite{lambert2008collision}, the integral (\ref{CollisonProbIndicator}) was approximated via empirical average using the state samples. For our computational framework, in addition to the state samples, the transient bivariate marginals $\xi^{\text{ego}}, \xi^{\text{non-ego}}$ are numerically available, which allow us to approximate (\ref{CollisonProbIndicator}) as a weighted sum. With the simulation set up described above, Fig. \ref{FigCollision} shows $p_{\text{collision}}(t)$ for $t\in[t_{0},t_{f}]$.

\subsection{Case Study: Barycentric Planning in the Belief Space}\label{subsecCaseStudy}

\begin{figure}[htpb]
\centering
 \includegraphics[width=0.65\linewidth]{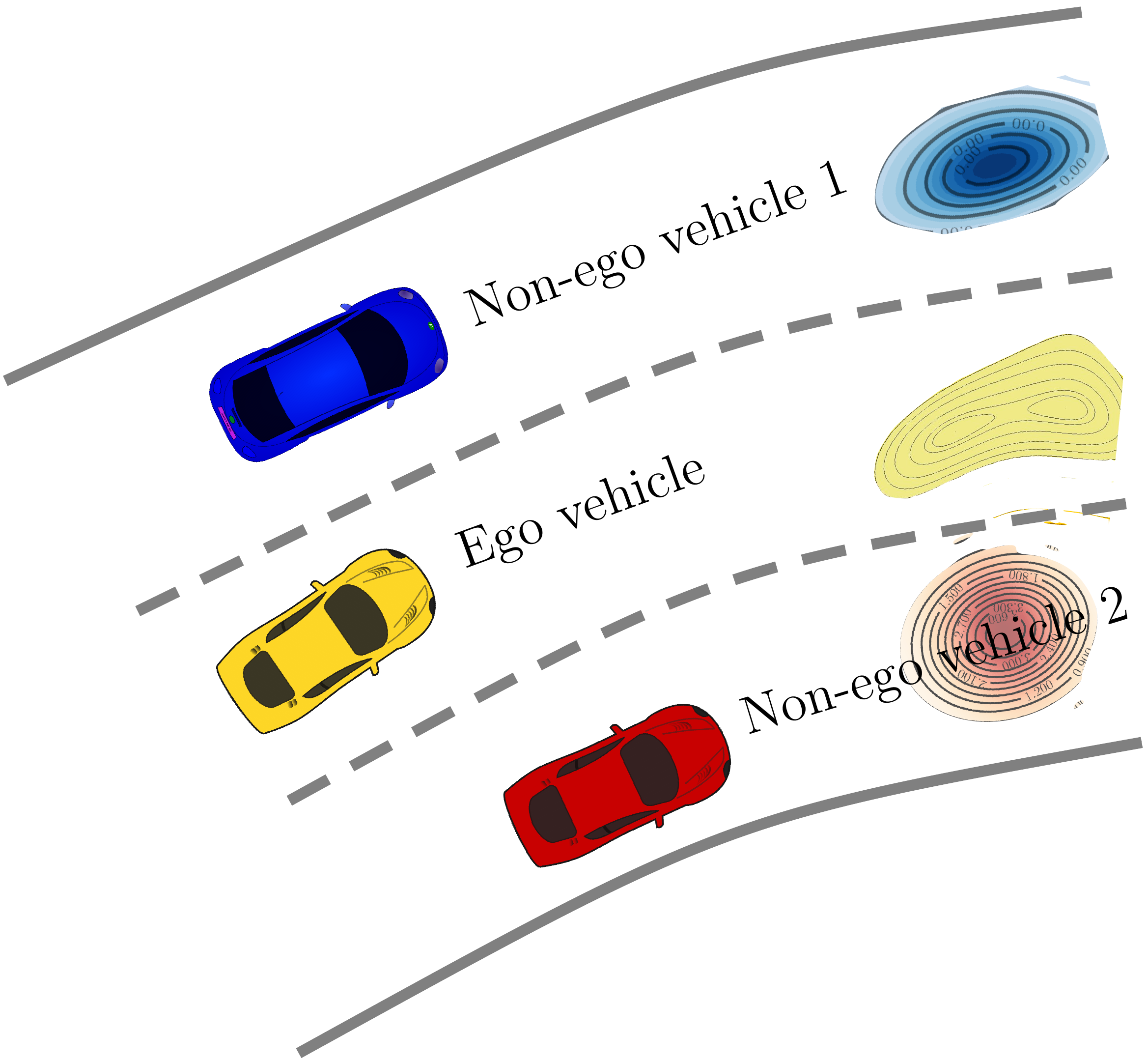}
         \caption{{\small{A schematic of the three lane highway driving scenario described in Section \ref{subsecCaseStudy}. As shown, the ego vehicle is in the middle lane while the two non-ego vehicles are in its left and right lanes, respectively. The three contour plots show the respective bivariate marginals $\xi^{\text{ego}}$, $\xi^{\text{non-ego (1)}}$, and $\xi^{\text{non-ego (2)}}$ in $(s,e_y)$ variables at an instance $t$ within the prediction horizon $[t_0,t_f]$.}}}
\vspace*{-0.14in}
\label{ThreeCarsThreeLanes}       
\end{figure}
To illustrate the scope of the proposed computational framework, we now consider a three lane (unidirectional) highway driving scenario as shown in Fig. \ref{ThreeCarsThreeLanes}. The ego vehicle, shown in the middle lane in Fig. \ref{ThreeCarsThreeLanes}, at time $t_{0}$, computes the evolution of its own joint $\rho^{\text{ego}}(\bm{x},t)$ and bivariate marginal $\xi^{\text{ego}}\left(s^{\text{ego}},e_{y}^{\text{ego}},t\right)$, as well as its prediction of the same for its neighboring non-ego vehicles, i.e., 
\[\rho^{\text{non-ego $(i)$}}(\bm{x},t), \: \xi^{\text{non-ego $(i)$}}\left(s^{\text{non-ego $(i)$}},e_{y}^{\text{non-ego $(i)$}},t\right), \: i\in\{1,2\},\] over a prediction horizon $[t_{0},t_{f}]$. Then,  the ideas in Section \ref{SimulationDynamic} allow the ego to predict the time-varying collision probabilities between each ego-non-ego pair. Moreover, if the ego would like to continue in its (here, middle) lane, then a natural question is how should it plan its probabilistically safest path in the $(s,e_y)$ variables over the horizon $[t_0,t_f]$ taking into account the stochastic predictions $\xi^{\text{non-ego $(i)$}}\left(s^{\text{non-ego $(i)$}},e_{y}^{\text{non-ego $(i)$}},t\right)$, $i\in\{1,2\}$, $t\in[t_0,t_f]$. Specifically, supposing that the non-egos continue in their respective lanes over the prediction horizon, the ego would like to plan a path in the $(s,e_y)$ marginal PDF or belief space, that is maximally separated from both the predicted belief trajectories $\xi^{\text{non-ego $(i)$}}\left(s^{\text{non-ego $(i)$}},e_{y}^{\text{non-ego $(i)$}},t\right)$, $i\in\{1,2\}$, thus minimizing the chance of collision with any of the non-ego vehicles for all $t\in[t_0,t_f]$.

\begin{figure*}[htpb]
\centering
 \includegraphics[width=\linewidth]{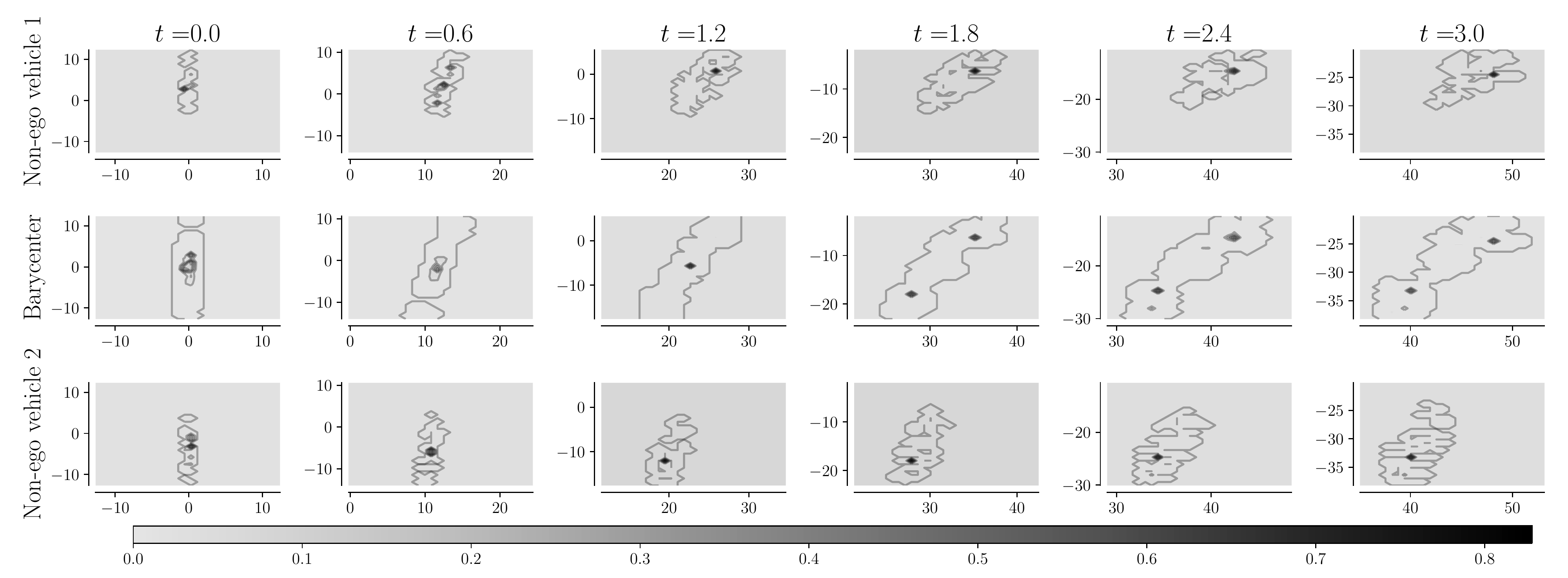}
         \caption{{\small{For the simulation set up detailed in Section \ref{subsecCaseStudy}, the evolution of the transient bivariate marginal state PDFs in $\left(s,e_y\right)$ variables for the non-ego vehicle 1 (top row), the non-ego vehicle 2 (bottom row), and their barycenters (middle row). The colorbar (\emph{light hue = small}, \emph{dark hue = large}) shows the values of these bivariate marginals.}}}
\vspace*{-0.2in}
\label{ThreeVehicleDynamicJoint}       
\end{figure*}

For a given bivariate marginal PDF $\xi(s,e_y)$, and for $i\in\{1,2\}$, let
\begin{align*}
W_{i}(t) \!:=\! W\!\!\left(\!\xi(s,e_y), \xi^{\text{non-ego $(i)$}}\!\left(\!s^{\text{non-ego $(i)$}},e_{y}^{\text{non-ego $(i)$}},t\right)\!\right)
\end{align*}
denote the 2-Wasserstein distance between $\xi$ and $\xi^{\text{non-ego $(i)$}}$ at time $t$. The 2-Wasserstein distance $W$ is a metric in the space of PDFs (or probability measures, in general), and it quantifies the minimal cost of transporting or reshaping one PDF to another; see e.g., \cite[Ch. 7]{villani2003topics}. For any two probability measures $\varsigma_{0},\varsigma_{1}$ with finite second moments, each supported on $\mathcal{M}\subseteq\mathbb{R}^{d}$, it is defined as
\begin{align}
W(\varsigma_{0},\varsigma_{1}) \!:= \!\left[\underset{\wp\in\Pi(\varsigma_{0},\varsigma_{1})}{\inf}\int_{\mathcal{M}\times\mathcal{M}}\!\!\|\bm{x} - \bm{y}\|_{2}^{2} \:\differential\wp(\bm{x},\bm{y}) \right]^{\frac{1}{2}},
\label{defW}	
\end{align}
where $\Pi(\varsigma_{0},\varsigma_{1})$ is the set of joint probability measures on $\mathcal{M}\times\mathcal{M}$ with marginals $\varsigma_{0},\varsigma_{1}$. Then, maximizing the distance from the non-ego marginal beliefs, as motivated above, amounts to computing the barycentric \cite{agueh2011barycenters} trajectory for $t\in[t_0,t_f]$ in the belief space: 
\begin{align}
\xi^{\text{bary}}(s,e_y,t) := \underset{\xi\in\mathcal{P}_{2}\left(\mathbb{R}^2\right)}{\arg\inf}\displaystyle\sum_{i=1}^{2}\lambda_{i}(t)\left(W_{i}(t)\right)^{2},
\label{WassBary}	
\end{align}
where the weights $\lambda_{i}(t)$ codify the ego vehicle's confidence in its stochastic prediction $\xi^{\text{non-ego $(i)$}}\left(s^{\text{non-ego $(i)$}},e_{y}^{\text{non-ego $(i)$}},t\right)$, and satisfy $\lambda_{i}(t)\geq 0, \sum_{i=1}^{2}\lambda_{i}(t)=1$. In (\ref{WassBary}), $\mathcal{P}_{2}\left(\mathbb{R}^2\right)$ is the set of all joint PDFs over $\mathbb{R}^{2}$ with finite second moments. Notice that (\ref{WassBary}) generalizes naturally to the case when each neighboring lane of ego has more than one non-ego vehicle. 

Intuitively, the barycenter $\xi^{\text{bary}}(s,e_y,t)$ is the weighted average of the beliefs $\bigg\{\xi^{\text{non-ego $(i)$}}\left(s^{\text{non-ego $(i)$}},e_{y}^{\text{non-ego $(i)$}},t\right)\bigg\}_{i=1}^{2}$. In the absence of uncertainties, the beliefs $\xi^{\text{non-ego $(i)$}}$ would be Dirac deltas, and the barycentric PDF trajectory would reduce to a Dirac delta trajectory supported on the graph of the weighted average of the non-ego state trajectories.

To numerically simulate the scenario shown in Fig. \ref{ThreeCarsThreeLanes}, let
\begin{align*}
\rho_{0}^{\text{ego}}(\bm{x}) &= \mathcal{N}\left(\bm{\mu}_{0}^{\text{ego}},\bm{\Sigma}_{0}^{\text{ego}}\right),\\
\rho_{0}^{\text{non-ego $(i)$}}(\bm{x}) &= \mathcal{N}\left(\bm{\mu}_{0}^{\text{non-ego $(i)$}},\bm{\Sigma}_{0}^{\text{non-ego $(i)$}}\right), \quad i\in\{1,2\},	
\end{align*}
with
\vspace*{-0.1in}
\begin{align*}
\bm{\mu}_{0}^{\text{ego}} &= \left(20, 0, 0, 0, 0, 0\right)^{\top}, \\
\bm{\Sigma}_{0}^{\text{ego}} &= {\rm{diag}}\left(1.11\times 10^{-3}, 1.11\times 10^{-3}, 1.23\times 10^{-8},\right. \\
&\left.\qquad\quad\;\; 2.78\times 10^{-6}, 1.11\times 10^{-3}, 1.11\times 10^{-1}\right),\\
\bm{\mu}_{0}^{\text{non-ego $(i)$}} &= \left(20, 0, 0, 0, (-1)^{i+1}\times 3.7, 0\right)^{\top},\\
\bm{\Sigma}_{0}^{\text{non-ego $(1)$}} &= \bm{\Sigma}_{0}^{\text{non-ego $(2)$}} = {\rm{diag}}\left(0.11, 0.11, 1.23\times10^{-6},\right.\\
&\left.\qquad\qquad\qquad\qquad\quad 2.5\times10^{-5}, 11.11, 0.11\right).	
\end{align*}
Following Section \ref{FeedbackSynthesis}, we compute the $v_{x}$ velocity-hold trim conditions for each of the three vehicles, and perform the offline linearized MPC synthesis about the same as in Section \ref{SimulationDynamic}. We use $N=200$ random samples for each of the three initial joint state PDFs given above, and as in Section \ref{SimulationDynamic}, generate weighted scattered point clouds $\{\bm{x}^{\text{ego},i}(t),\rho^{\text{ego},i}(t)\}_{i=1}^{N}$, $\{\bm{x}^{\text{non-ego 1},i}(t),\rho^{\text{non-ego 1},i}(t)\}_{i=1}^{N}$, $\{\bm{x}^{\text{non-ego 2},i}(t),\rho^{\text{non-ego 2},i}(t)\}_{i=1}^{N}$, which are then propagated via the Liouville PDE in Section \ref{SecStocReachability} over $t\in[t_{0},t_{f}] \equiv [0,3]$ with the respective explicit MPC feedbacks in the loop. As before, we then use the resulting joint PDFs to compute the time-varying bivariate marginals $\xi^{\text{ego}}, \xi^{\text{non-ego 1}}, \xi^{\text{non-ego 2}}$. To compute the barycentric marginals in (\ref{WassBary}), we employ the multi-marginal entropy-regularized optimal transport proposed in \cite[Section 4.2]{benamou2015iterative} with $\lambda_1(t)=\lambda_2(t) = 1/2$ for all $t\in[t_{0},t_{f}]$.

\begin{figure*}[tpb]
\centering
 \includegraphics[width=0.82\linewidth]{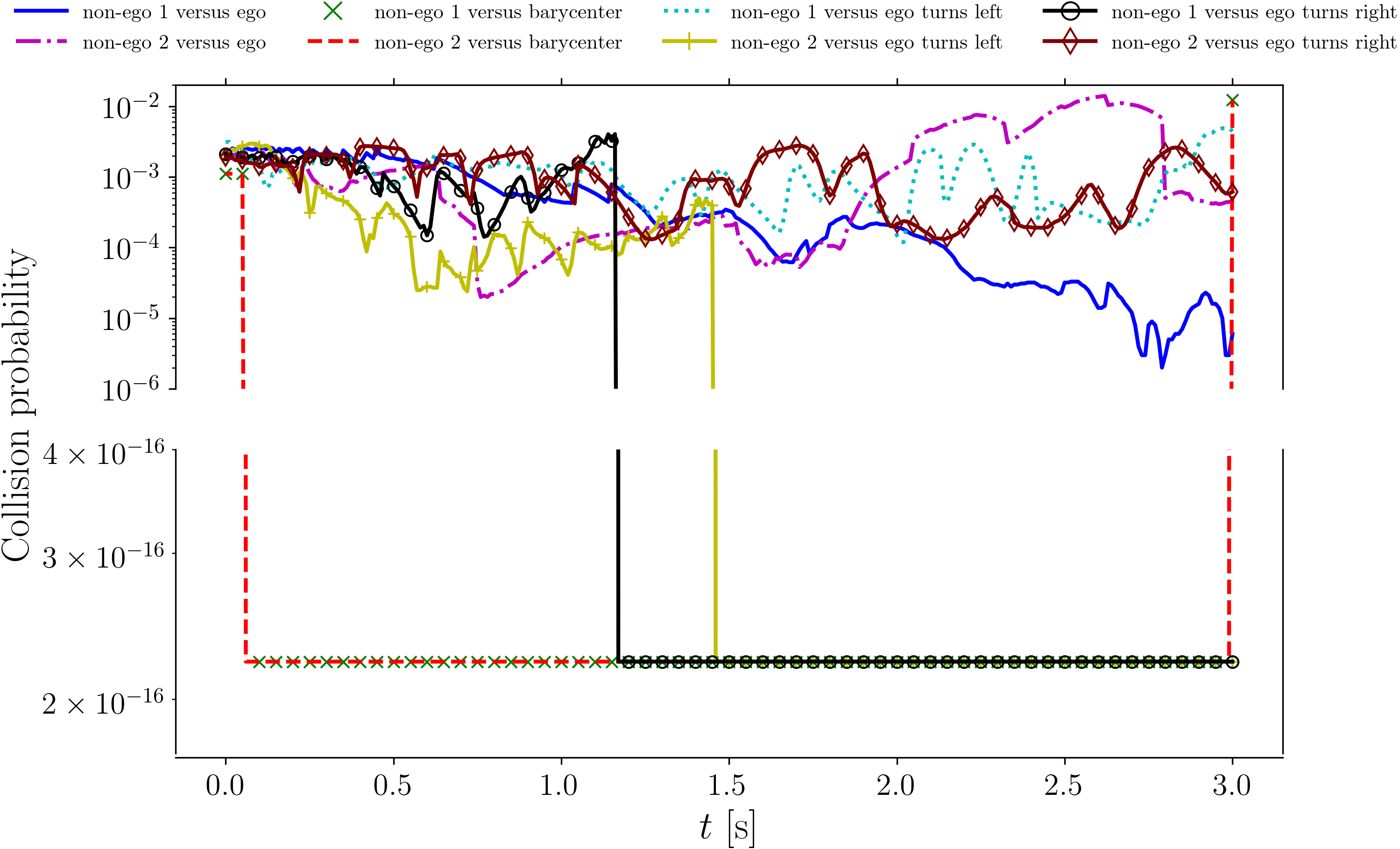}
         \caption{{\small{The time-varying collision probabilities for the simulation set up detailed in Section \ref{subsecCaseStudy}, plotted with a broken vertical axis as some collision probabilities drop to (machine precision) zero.}}}
\vspace*{-0.18in}
\label{FigCollisionProbThreeVehicle}       
\end{figure*}

Fig. \ref{ThreeVehicleDynamicJoint} shows the evolution of the transient marginal state PDFs in $(s,e_y)$ variables for the non-ego vehicles, as well as their barycenters. Fig. \ref{FigCollisionProbThreeVehicle} highlights that the collision probabilities between the barycenter-non-ego pairs are negligible compared to the collision probabilities between the ego-non-ego pairs, as expected. Also shown in Fig. \ref{FigCollisionProbThreeVehicle} are the ego versus non-ego collision probabilities when the ego (starting from stochastic initial conditions in the middle lane, as before) turns either to its left or to its right lane.

\begin{figure}[htpb]
\centering
 \includegraphics[width=\linewidth]{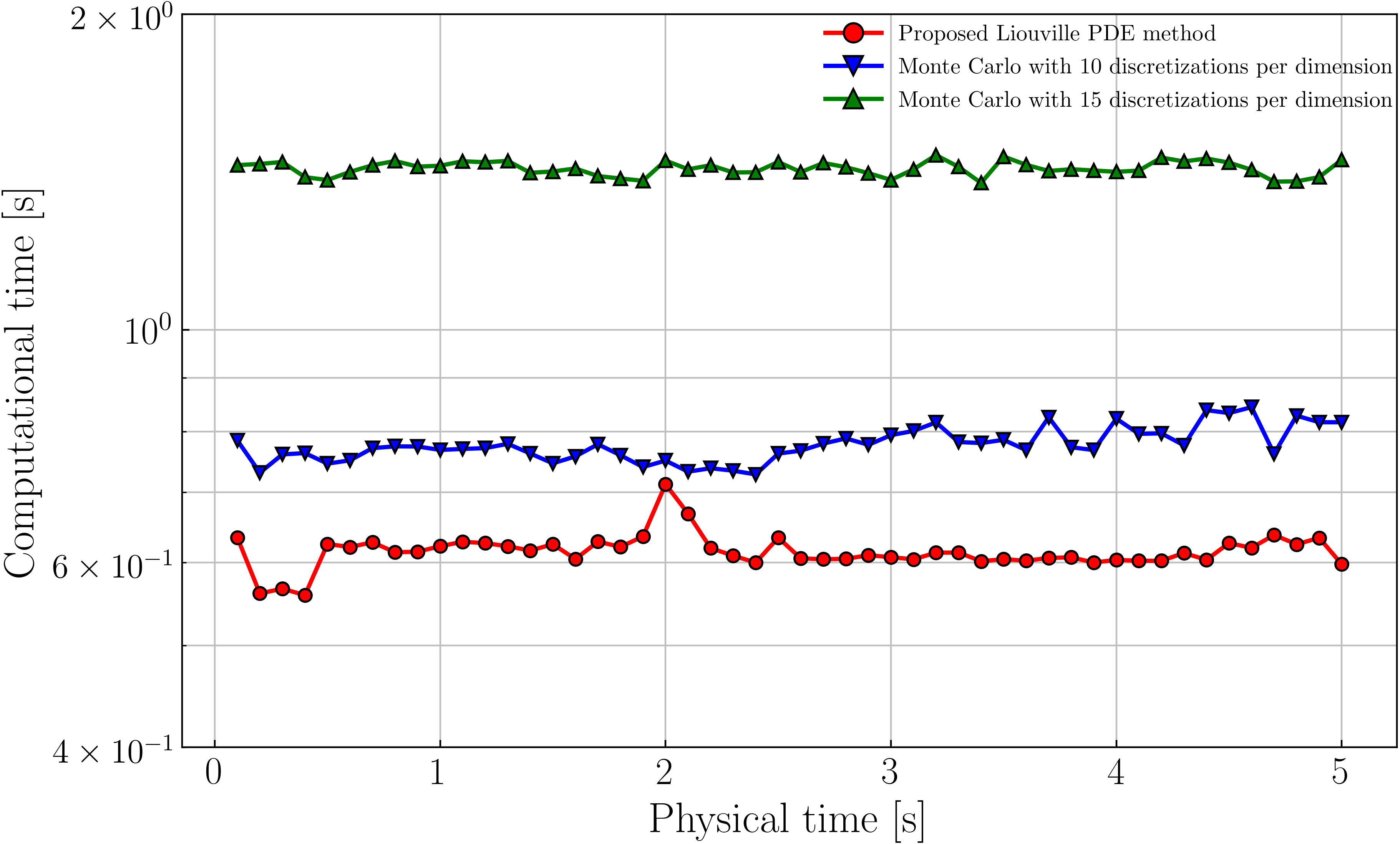}
         \caption{{\small{Comparison of the computational times for calculating the transient joint PDFs of the ego vehicle in Section \ref{SimulationKinematic} using the proposed method and the standard Monte Carlo.}}}
\vspace*{-0.2in}
\label{CompTimeLiouvilleVersusMC}       
\end{figure}

\begin{figure}[htpb]
\centering
 \includegraphics[width=\linewidth]{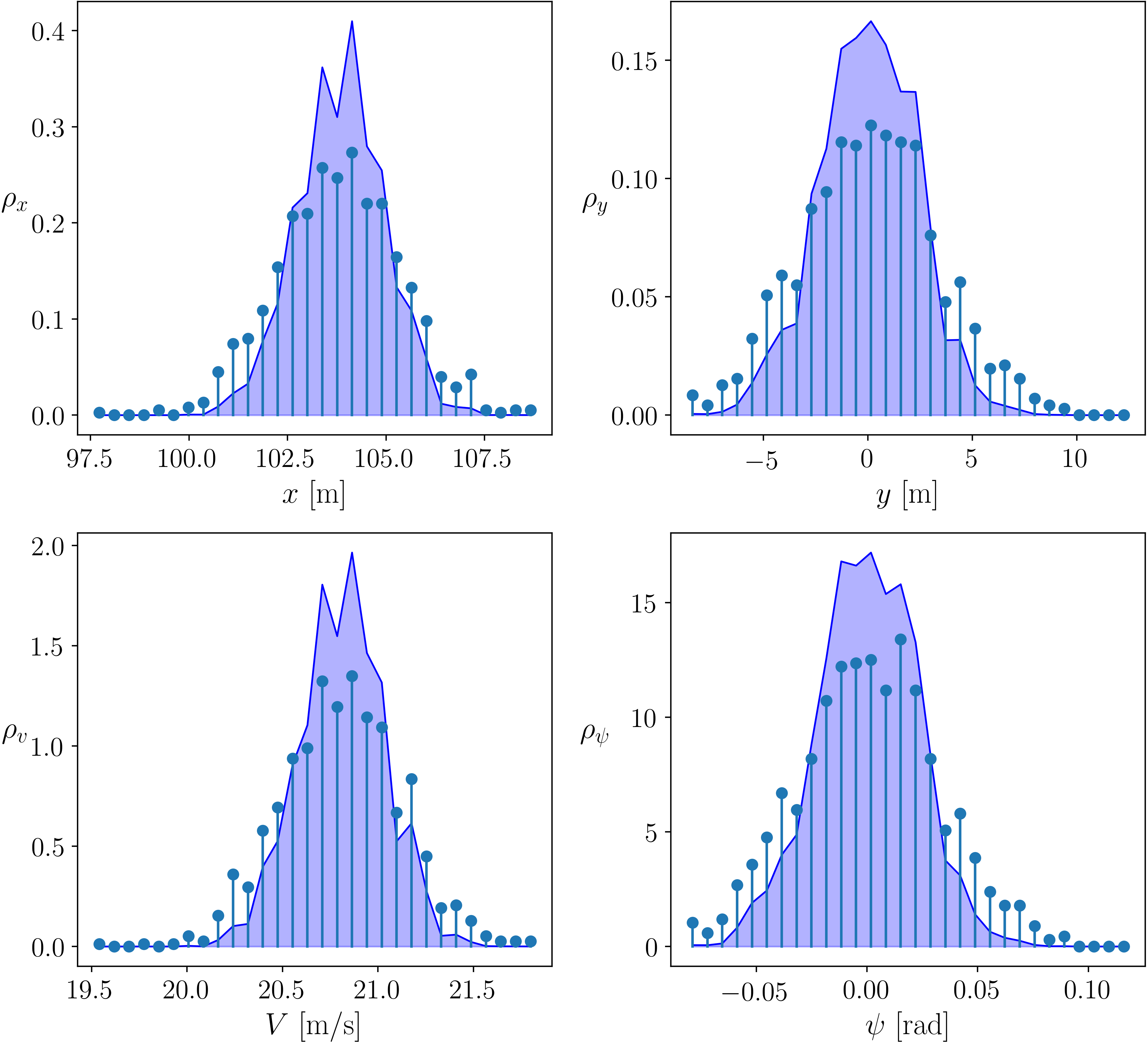}
         \caption{{\small{Comparison of the univariate marginal PDFs of the ego vehicle in Section \ref{SimulationKinematic} at $t=5$ s, using the proposed method (\emph{solid lines with filled area}) and the standard Monte Carlo (\emph{stem plot}).}}}
\vspace*{-0.2in}
\label{AccuracyLiouvilleVersusMC}       
\end{figure}

\subsection{Comparison with the Monte Carlo Method}\label{subsecComparisonMC}
We now provide a quantitative comparison of the proposed Liouville PDE method with the standard Monte Carlo. To this end, we compare the computational time and accuracy for the ego vehicle's joint PDF prediction via these two methods in context of the simulation set up in Section \ref{SimulationKinematic}.

Specifically, Fig. \ref{CompTimeLiouvilleVersusMC} shows that the PDE characteristic-based PDF propagation, as outlined in Section \ref{SecStocReachability}, is computationally faster compared to the standard Monte Carlo method. The latter requires constructing a grid over the (in this case, four dimensional) state space to approximate the joint state PDF. Even though the state samples can be propagated in a meshless manner, grids are necessary in the Monte Carlo method for constructing histograms representing the piecewise constant function approximation for the instantaneous joint PDFs. As mentioned in Section \ref{SecStocReachability}, the proposed method--not being a function approximation--obviates the construction of grid, but requires integrating an extra ODE per sample; see (\ref{CharODE}). 

In Fig. \ref{CompTimeLiouvilleVersusMC}, we used $N=1000$ samples with the same parameters and initial joint PDFs detailed in Section \ref{SimulationKinematic}. For the same data, two Monte Carlo simulations were performed: one with 10, and another with 15 uniform discretizations per state dimension, between the transient minimum and maximum along each state coordinate. This resulted in histograms supported over $10^{4}$ and $15^{4}$ orthotopes, respectively. Fig. \ref{CompTimeLiouvilleVersusMC} also shows that finer discretization for Monte Carlo requires more computational time, as expected.

Fig. \ref{AccuracyLiouvilleVersusMC} compares the univariate marginal PDFs for the ego vehicle at $t=5$ computed via the proposed method and the standard Monte Carlo. While the overall trends are similar, we point out that the solid lines in Fig. \ref{AccuracyLiouvilleVersusMC} are obtained by numerically integrating the joint PDF values via (\ref{CharODE}), which are \emph{exact modulo the floating point arithmetic}. In contrast, the Monte Carlo marginals are obtained by numerically integrating the \emph{piecewise constant approximations of the joint PDFs}, which are in turn obtained by constructing histograms over the propagated states. Of course, finer discretization in Monte Carlo can increase its accuracy but at the expense of computational time; cf. Fig. \ref{CompTimeLiouvilleVersusMC}.


\section{Conclusions}\label{secConclusions}
We have proposed a density-based stochastic reachability computation framework for occupancy prediction in automated driving. The main idea is to evolve the joint state PDFs along the characteristic curves of the underlying transport PDE via gridless computation. These transient joint PDFs are then used to compute the marginals, and the time-varying collision probabilities for the closed-loop vehicle dynamics. We envision that such computational results can be passed to a higher level decision-making module such as driver-assistance or early collision warning systems. Although our simulations focused on the multi-lane highway driving scenarios, the overall framework can be adapted to applications that involve interactions with other vehicles and pedestrians in complex traffic environment. Another possible direction of future research is to incorporate the process noise (see e.g., \cite{caluya2019proximal,caluya2019gradient}) in the closed-loop vehicle dynamics.

\end{document}